\documentclass[leqno]{article}
\usepackage{amsmath}
\usepackage{amssymb,amscd}
\usepackage{verbatim}

\renewcommand{\baselinestretch}{1.2}
\begin{document}
\hyphenation{semi-per-fect Grothen-dieck}
\newtheorem{Lemma}{Lemma}[section]
\newtheorem{Th}[Lemma]{Theorem}
\newtheorem{Prop}[Lemma]{Proposition}
\newtheorem{OP}[Lemma]{Open Problem}
\newtheorem{Cor}[Lemma]{Corollary}
\newtheorem{Fact}[Lemma]{Fact}
\newtheorem{Def}[Lemma]{Definition}
\newtheorem{Not}[Lemma]{Notation}
\newtheorem{Ex}[Lemma]{Example}
\newtheorem{Exs}[Lemma]{Examples}
\newtheorem{Rem}[Lemma]{Remark}
\newenvironment{Remarks}{\noindent {\bf Remarks.}\ }{}
\newtheorem{Remark}[Lemma]{Remark}
\newenvironment{Proof}{{\sc Proof.}\ }{~\rule{1ex}{1ex}\vspace{0.5truecm}}
\newenvironment{ProofTh}{{\sc Proof of the Theorem.}\ }{~\rule{1ex}{1ex}\vspace{0.5truecm}}
\newenvironment{Proofmain}{\noindent{\textbf{ Proof  Theorem \ref{main}.}}\ }{~\rule{1ex}{1ex}\vspace{0.5truecm}}
\newcommand{\End}{\mbox{\rm End}}
\newcommand{\K}{\mbox{\rm K.dim}}
\newcommand{\Hom}{\mbox{\rm Hom}}
\newcommand{\infsupp}{\mathrm{inf\mbox{-}supp}\,}
\newcommand{\finsupp}{\mathrm{fin\mbox{-}supp}\,}
\newcommand{\Ext}{\mbox{\rm Ext}}
\newcommand{\supp}{\mbox{\rm supp}\,}
\newcommand{\Supp}{\mbox{\rm Supp}\,}
\newcommand{\Max}{\mbox{\rm Max}}
\newcommand{\cl}{\mbox{\rm cl}}
\newcommand{\add}{\mbox{\rm add}}
\newcommand{\Inv}{\mbox{\rm Inv}}
\newcommand{\rk}{\mbox{\rm rk}}
\newcommand{\Tr}{\mbox{\rm Tr}}
\newcommand{\Sat}{\mbox{\bf Sat}}
\newcommand{\card}{\mbox{\rm card}}
\newcommand{\Ann}{\mbox{\rm Ann}}
\newcommand{\proj}{\mbox{\rm proj-}}
\newcommand{\codim}{\mbox{\rm codim}}
\newcommand{\B}{\mathcal{B}}
\newcommand{\Scal}{\mathcal{S}}
\newcommand{\Cong}{\mbox{\rm Cong}}
\newcommand{\Spec}{\mbox{\rm Spec-}}
\newcommand{\coker}{\mbox{\rm coker}}
\newcommand{\Cl}{\mbox{\rm Cl}}
\newcommand{\Ses}{\mbox{\rm Ses}}
\newcommand{\im}{\mbox{\rm Im}}
\newcommand{\Cal}[1]{{\cal #1}}
\newcommand{\+}{\oplus}
\newcommand{\N}{\mathbb N}
\newcommand{\No}{{\mathbb N}_0}
\newcommand{\Z}{\mathbb{Z}}
\newcommand{\Q}{\mathbb{Q}}
\newcommand{\C}{\mathbb{C}}
\newcommand{\T}{\mathbb{T}}
\newcommand{\R}{\mathbb{R}}
\newcommand{\notsim}{\sim\makebox[0\width][r]{$\slash\;\,$}}
\newcommand{\Mod}{\mbox{\rm Mod-}}
\newcommand{\lmod}{\mbox{\rm -mod}}
\newcommand{\mspec}{\mbox{\rm Max}}
\renewcommand{\dim}{\mathrm{dim}\,}

\title {Big projective modules over noetherian semilocal rings}

\author{Dolors Herbera\thanks{Partially supported by MEC-DGESIC
(Spain) through Project MTM2005-00934, and by the Comissionat Per
Universitats i Recerca de la Generalitat de Catalunya through
Project 2005SGR00206.  \protect\newline This paper was started
when both authors were working within the Research Programme
 on Discrete and Continuous methods of Ring Theory at the CRM, Barcelona
 (Spain). They thank their host for its hospitality.}\\
 Departament de Matem\`atiques, \\
Universitat Aut\`onoma de Barcelona, \\ 08193 Bellaterra
(Barcelona), Spain\\ e-mail: dolors@mat.uab.cat  \and Pavel P\v
r\'\i hoda\thanks{At CRM supported by grant SB 2005-0182 from MEC, recently supported by research project MSM 0021620839.} \\ Charles University, Faculty of Mathematics and Physic, \\Department of Algebra, Sokolovsk\'a 83,\\ 186 75 Praha 8,
Czech Republic\\ e-mail: prihoda@karlin.mff.cuni.cz}

\date{\phantom{ciao}}

\maketitle

\begin{abstract} We prove that for a noetherian semilocal ring $R$
with exactly $k$ isomorphism classes of simple right  modules the
monoid $V^*(R)$ of isomorphism classes of countably generated
projective right (left) modules, viewed as a submonoid of
$V^*(R/J(R))$, is isomorphic to the monoid of solutions in $(\No
\cup\{\infty\})^k$ of a system consisting of congruences and
diophantine linear equations. The converse also holds, that is, if
$M$ is a submonoid of $(\No \cup\{\infty\})^k$ containing an order
unit $(n_1,\dots ,n_k)$ of $\No^k$ which is the set of  solutions of
a system of congruences and linear diophantine equations then it can
be realized as $V^*(R)$ for a noetherian semilocal ring such that
$R/J(R)\cong M_{n_1}(D_1)\times \cdots \times M_{n_k}(D_k)$ for
suitable division rings $D_1,\dots, D_k$.
\end{abstract}

A theorem of Kaplansky states that, for any ring $R$, a projective
right $R$-module is a direct sum of countably generated projective
right $R$-modules. This reduces the study of direct summands
of $R^{(I)}$, where $I$ denotes an arbitrary set, to
the study of direct sum decomposition of $R^{(\omega)}$ or,
equivalently, to the study of countably generated projective right
$R$-modules.

The commutative monoid $V(R)$ of isomorphism classes of finitely
generated projective right $R$-modules, with the addition induced by
the direct sum of modules, encodes  the direct sum behavior of
finite direct sums of  finitely generated projective right
$R$-modules. Similarly, the monoid $V^*(R)$ of isomorphism classes
of countably generated projective right $R$-modules, with the
addition induced by the direct sum of modules, encodes the
direct-sum behavior of countably generated projective modules. In
this paper we characterize the monoids  that can be realized as
$V^*(R)$ for $R$  a noetherian semilocal ring.

A ring $R$ is said to be semilocal if it is semisimple artinian
modulo its Jacobson radical $J(R)$. To fix notation, we assume that
$R/J(R)\cong M_{n_1}(D_1)\times \cdots \times M_{n_k}(D_k)$ where
$D_1,\dots ,D_k$ are division rings uniquely determined up to
isomorphism. The monoids $V(R)$ and $V^*(R)$ can be viewed as
submonoids of $V(R/J(R))\cong \No^k$ and $V^*(R/J(R))\cong (\No\cup
\{\infty\})^k=(\No^*)^k$, respectively; the class of $R$ corresponds
to the element $(n_1,\dots ,n_k)\in \No^k$. The submonoids of
$\No^k$ containing $(n_1,\dots ,n_k)$ that can be realized as $V(R)$
for a semilocal ring $R$ were characterized in \cite{FH} as the set
of solutions in $\No^k$ of systems of diophantine equations of the
form
\[D\left(\begin{array}{c}t_1\\\vdots \\ t_k\end{array}\right)\in \left(\begin{array}{c}m_1\No^*\\
\vdots \\ m_n\No^*\end{array} \right)\quad \qquad \mbox{and} \qquad
E_1\left(\begin{array}{c}t_1\\\vdots \\
t_k\end{array}\right)=E_2\left(\begin{array}{c}t_1\\\vdots \\
t_k\end{array}\right)\quad (1)\]  where the  coefficients of the
matrices $D$, $E_1$ and $E_2$ as well as $m_1,\dots ,m_n$ are
elements of $\No$. Such submonoids of $\No^k$ are called  full
affine submonoids (cf. Definition~\ref{deffull} and
Proposition~\ref{relace}). This terminology was introduced by
Hochster in \cite{Hochster}. However full affine monoids appear in
different contexts with different names. In the setting of
commutative noetherian rings they are also called positive normal
monoids, see \cite{bruns2}. Such monoids also appear in
generalizations of the multiplicative ideal theory where they are
called finitely generated Krull monoids see, for example,
\cite{chuinard}.

In this paper we show that the submonoids of $(\No^*)^k$ that can be
realized as $V^*(R)$ for a noetherian semilocal ring $R$ are
precisely the sets of solutions in $(\No^*)^k$ of systems of  type
$(1)$. We refer the reader to Theorem~\ref{main} for the precise
statement. Most of the paper is devoted to the proof of
Theorem~\ref{main} which has, essentially, two quite different
parts. A more ring theoretical one, in which we provide the
necessary tools to construct noetherian semilocal rings with
prescribed monoid $V^*(R)$. Our key idea is to use a well known
theorem due to Milnor stating that, under mild restrictions, the
category of right projective modules over a pullback of rings is
equivalent to the pullback of the categories of projective modules.
Surprisingly enough, just considering pullbacks of semilocal
principal ideal  domains (or just noetherian semilocal rings such
that all projective modules are free) and semisimple artinian rings
a rich supply of noetherian semilocal rings $R$ with non-trivial
$V^*(R)$ is obtained.

The second part of the paper (and of the proof of
Theorem~\ref{main}) deals with submonoids of $(\No^*)^k$. Our
starting point are the results in \cite{fair} where it was proven
that, for a noetherian semilocal ring, $V^*(R)$ is built up from a
finite collection of full affine submonoids of $\No ^{r_1},\dots,
\No^{r_m}$, respectively, where $r_i\le k$, chosen in a
\emph{compatible} way. These monoids  are placed in the finite
supports of the elements of $V^*(R)$ or, better saying, in the
complementary of the infinite supports of the elements of $V^*(R)$,
see Definition~\ref{supports} for the unexplained terminology. In
the paper, we make an abstraction of this type of monoid  by
introducing the concept of \emph{full affine system of supports} in
Definition~\ref{defsystems}, then $V^*(R)$, viewed as a submonoid of
$(\No^*)^k$, is given by a \emph{full affine system of supports}.
Our main result in this part of the paper shows that the monoids
given by a full affine system of supports are precisely the
solutions in $(\No^*)^k$ of systems of the form $(1)$.

We stress the fact that though the description of these submonoids
of $(\No^*)^k$ as  sets of solutions of a system of equations is
very elegant, and it extends nicely the characterization for the
case of finitely generated projective modules, the one given by the
systems of supports seems to give  a better idea of the complexity
of the monoids we are working with. The contrast with the
commutative situation is quite striking,    all projective modules
over a commutative semilocal indecomposable ring are free \cite{Hi}.
On the other hand, as we explain along the paper and especially in
\S \ref{final} (cf. Definition~\ref{cuntzmonoid} and
Examples~\ref{noniso}), the noetherian situation is simpler than the
general one.

Our interest on semilocal rings stems from the fact that many
classes of \emph{small modules} have a semilocal endomorphism ring.
For example, artinian modules or, more generally, modules with
finite Goldie and dual Goldie dimension \cite{HS}, finitely
presented modules over a local ring or, more generally, finitely
presented modules over a semilocal ring  are classes of modules with
a semilocal endomorphism ring \cite{limit}. We refer to the
monograph \cite{libro} as a source to read  about, the good and the
not so good, properties of modules with a semilocal endomorphism
ring.

A description of the projective modules over the endomorphism ring
is a first step towards understanding (part of) the category
$\mathrm{Add}\, (M)$ of direct summands of any direct sum of copies
of $M$. Wiegand in \cite{W} proved that all monoids of solutions in
$\No^k$ of systems of the form $(2)$ can be realized as $V(R)$ when
$R$ is the endomorphism ring of a finitely presented module over a
noetherian semilocal ring, or the endomorphism ring of an artinian
module. Yakovlev \cite{yakovlev1, yakovlev2} proved the same kind of
result for semilocal endomorphisms rings of certain classes of
torsion free abelian groups of finite rank. For further information,
see the   survey  paper \cite{WW}. Our results give a new twist to
the situation, as they  indicate that when considering countable
direct sums of such modules a rich supply of new direct summands
might appear.

Let us mention a connection between projective modules over
noetherian semilocal ring and  integral representation theory. In
\cite{BCK} the so called generalized lattices were investigated. For
a Dedekind domain $D$ with a quotient field $K$ we consider an order
$R$ in a separable $K$-algebra. An $R$-module $M$ is said to be a
generalized $R$-lattice provided it is projective as a $D$-module.
If $M$ is also finitely generated,  $M$ is a lattice over $R$. $R$
is said to be of finite lattice  type if there exist only finitely
many indecomposable lattices up to isomorphism. Suppose that $R$ is
of finite lattice type and let $A$ be the direct sum of a
representative set of isomorphism classes of indecomposable
lattices. By \cite{BCK}, the category of generalized lattices over
$R$ and the category of projective modules over ${\rm End}_R(A)$ are
equivalent. For any maximal ideal $\mathcal{M}$ of $D$ let
$R_{(\mathcal{M})}$ be the localization of $R$ in $D \setminus
\mathcal{M}$ and let $R_{(0)}$ be the localization of $R$ in $D
\setminus \{0\}$. Then monoid homomorphisms $$V^*({\rm End}_R(A))
\to V^{*}({\rm End}_{R_{(\mathcal{M})}}(A \otimes_{R}
R_{(\mathcal{M})})) \to V^*({\rm End}_{R_{(0)}} (A \otimes_{R}
R_{(0)}))$$ give  approximations of generalized lattices over $R$
by projective modules over the  noetherian semilocal rings ${\rm
End}_{R_{(\mathcal{M})}}(A \otimes_{R} R_{(\mathcal{M})})$ and an artinian ring ${\rm
End}_{R_{(0)}} (A \otimes_{R} R_{(0)})$. For further results on generalized
lattices see \cite{Rump}.

The paper is structured as follows, in \S \ref{monoids} we introduce
the basic language used throughout the paper. We  describe the
monoids of projective modules,  specializing to a semisimple
artinian ring, we recall the results needed to understand the
relation between these monoids when considered over $R$ and over
$R/J(R)$ emphasizing on  the particular case of semilocal rings. In
\S \ref{noetherian} we specialize to the  noetherian case; we recall
the results from \cite{fair}  essential for our investigation and we
state our main Theorem~\ref{main}. Sections \ref{PID},
\ref{pullbacks} and \ref{construction} deal with the realization
part of the proof of Theorem~\ref{main}; \S \ref{PID} shows how to
construct  principal ideal domains with prescribed semisimple factor
modulo the Jacobson radical, in \S \ref{pullbacks} we provide all
the results we need on ring pullbacks in order to be able  to
realized the monoids we want as $V^*(R)$ of semilocal noetherian
algebras in \S \ref{construction}.

In section \ref{solving} we turn towards  monoids. We recall some
basics on full affine monoids, and we prove the auxiliary results
that will allow us to conclude the proof of Theorem~\ref{main} in \S
\ref{systems}. In section \ref{final} we mostly present some
examples to illustrate consequences of Theorem~\ref{main} and to
show the potential of the pullback constructions also outside the
noetherian setting. For example, in \ref{noniso}, we construct a
semilocal ring such that all projective left $R$-modules are free
while $R$ has a nonzero (infinitely generated) right projective
module that is not a generator, so it is not a direct sum of
finitely generated projective modules. Such an example also shows
that the notion of p-connected ring is not left-right symmetric;
this answers in the negative a question in \cite[p. 310]{FS}. Recall
that, according to Bass \cite{bass}, a ring is (left) $p$-connected
if every nonzero left projective module is a generator

\bigskip

All rings have $1$, ring morphisms and modules are unital. We shall
usually consider right modules.

Our convention is $\N =\{1,2,\dots \}$, and we denote the
nonnegative integers by $\No =\{0, 1,2, \dots \}$.

Another basic object in this paper is the monoid
$(\mathbb{N}_0^*,+,0)$ whose underlying set is $\mathbb{N}_0 \cup
\{\infty\}$, the operation $+$ is the extension of addition of
non-negative integers by the rule $\infty+x = x + \infty = \infty$.
Sometimes we will be also interested in extending the product of
$\No$ to $\No^*$ by setting $\infty\cdot 0=0$ and $\infty \cdot
n=\infty$ for any $n\in \No^*\setminus \{0\}.$

For any right $R$-module $M$ the {\em trace} of $M$ in $R$ is the two
sided ideal of $R$
\[\mathrm{Tr}_R(M)=\mathrm{Tr}\, (M)=\sum_{f\in
\mathrm{Hom}_R(M,R)}f(M).\] If $X\subseteq M$ then  we denote the
{\em right annihilator} of $X$ by
\[r_R(X)=\{r\in R\mid mr=0\mbox{ for any } m\in X\}.\]
If $N$ is a left $R$-module and $Y\subseteq N$ then we denote the
{\em left annihilator} of $Y$ by
\[l_R(Y)=\{r\in R\mid rn=0\mbox{ for any } n\in Y\}.\]

\section{Monoids of projective modules} \label{monoids}

We recall the following definitions

\begin{Def} Let $(M,+,0)$ be a commutative additive monoid. An
element $x\in M$ is said to be an \emph{order unit} or \emph{an
archimedean element} of $M$ if for any $y\in M$ there exists $n\in
\N $ and $z\in M$ such that $nx=y+z$.

The monoid $M$ is said to be \emph{reduced} if, for any $x\in M$,
$x+y=0$ implies $x=0=y$.

Let $x,y\in M$. The relation $x\le y$ if and only if there exists
$z\in M$ such that $x+z=y$ is a preorder order on $M$ that is
called the \emph{algebraic order} or, more properly, the
\emph{algebraic preorder}.
\end{Def}

For example, any $x\in \No^*$  satisfies that $x\leq \infty$.  If
$k\ge 1$ the algebraic preorder of $\No ^k$ and in $(\No ^*)^k$ is
the component-wise order and it is a partial order.

Note that if  $M$ is a monoid preordered with the algebraic preorder
then all the elements must be positive.

\medskip

Let  $R$ be a ring.  We denote by $V(R)$   the monoid of isomorphism
classes of finitely generated projective right $R$-modules with the
operation induced by the direct sum. This is to say, if $P_1$ and
$P_2$ are finitely generated projective right $R$-modules then
$\langle P_1\rangle +\langle P_2\rangle=\langle P_1\oplus P_2\rangle
$. The monoid $V(R)$ is commutative, reduced and it has an order
unit $\langle R\rangle$. We usually think on $V(R)$ as a  monoid
preordered by the algebraic preorder.

Similarly, we define $V^* (R)$ to be the monoid of isomorphism
classes of countably generated projective right $R$-modules with the
sum induced by the direct sum. Clearly $V(R)$ is a preordered
submonoid of the preordered monoid $V^* (R)$.

The functor $\mathrm{Hom}_R(-,R)$ induces a monoid isomorphism
between $V(R)=V(R_R)$ and the monoid of isomorphism classes of
finitely generated left projective modules $V({}_RR)$. This is no
longer true for countably generated projective modules so, in
general, $V^*(R)=V^*(R_R)$ is not isomorphic to $V^*({}_RR)$, cf.
Example~\ref{noniso}.

 If $\varphi \colon R_1\to R_2$ is a ring morphism then the
functor $-\otimes _{R_1} R_2$ induces a morphism of monoids with
order unit $V(\varphi)\colon V(R_1)\to V(R_2)$ and a morphism of
monoids $V^*(\varphi) \colon V^* (R_1)\to V^* (R_2)$. Both morphisms
are given by the formula $\langle P\rangle\mapsto \langle P\otimes
_{R_1} R_2\rangle $. Another useful way to describe these monoid
morphisms is describing projective modules via idempotent matrices.

Let $P_{R_1}$ be a finitely generated projective right $R_1$
modules. There exist $n\in \N$ and an idempotent matrix $E\in
M_n(R_1)$ such that $P\cong ER_1^n$, then $V(\varphi)(\langle
ER_1^n\rangle )=\langle M_n(\varphi)(E)R_2^n\rangle $ where
$M_n(\varphi)\colon M_n(R_1)\to M_n(R_2)$ is the map defined by
$M_n(\varphi)(a_{ij})=(\varphi (a_{ij}))$. One proceeds similarly
with the countably generated projective right $R_1$-modules taking
instead of finite matrices elements in $\mathrm{CFM}(R_1)$ and
$\mathrm{CFM}(R_2)$  the rings of (countable) column finite matrices
with entries in $R_1$ and $R_2$ respectively.

\subsection{The semisimple artinian case} \label{semisimple} Let  $R$ be  a semisimple artinian  ring.
Therefore, there exist $k\in \N$, $n_1,\dots ,n_k\in \N$, $D_1,
\dots ,D_k$ division rings and  an isomorphism $\varphi\colon R\to
M_{n_1}(D_1)\times \dots \times M_{n_k}(D_k)$ with kernel $J(R)$.

Let  $(V_1, \dots ,V_k)$  be an ordered set of representatives of
the isomorphism classes of simple right $R$-modules such that $\End
_R(V_i)\cong D_i$ and, hence, $\mathrm{dim}\, ({}_{D_i}V_i)=n_i$ for
$i=1,\dots ,k$. If $P_R$ is a finitely generated projective module
then $P_R\cong V_1^{x_1}\oplus \dots \oplus V_k^{x_k}$. The
assignment $\langle P\rangle \mapsto (x_1,\dots ,x_k)\in \No ^k$
induces an isomorphism of monoids  $\dim _\varphi \colon V(R)\to \No
^k$. Since $\dim _\varphi (\langle R\rangle )=(n_1,\dots ,n_k)$,
taking $(n_1,\dots ,n_k)$ as the order unit of $\No ^k$, $\dim
_\varphi$ becomes an isomorphism of monoids with order unit. We call
$\dim _\varphi (\langle P\rangle )$ or, by abuse of notation $\dim
_\varphi (P)$, \emph{the dimension vector of the (finitely
generated) projective module $P$}.

The morphism $\dim _\varphi$ extends to a monoid morphism $\dim
_\varphi \colon V^* (R)\to (\No ^*) ^k$ by setting $\dim_{\varphi} (\langle
V_i^{(\omega)}\rangle )=(0,\dots,\infty ^{i)},\dots ,0)$ for
$i=1,\dots ,k$. Again, we call $\dim _\varphi(\langle P\rangle )$
\emph{the dimension vector of the (countably generated) projective
module $P$}.

Throughout the paper, especially in sections \ref{construction} and
\ref{final}, it is important to keep in mind how to compute
dimension vectors in terms of idempotent matrices. If $P$ is a
finitely generated (countably generated) right projective module
such that $\dim_\varphi (P)=(x_1,\dots ,x_k)$ then $P\cong
(E_1,\dots, E_k)\cdot F$ where $F$ is a finitely generated
(countably generated) free right $R$-module and $E_i$ are idempotent
matrices over $M_{n_i}(D_i)$ (over $\mathrm{CFM} (D_i)$) such that
$\mathrm{rank}_{D_i}(E_i)=x_i$ for $i=1,\dots,k$.

Finally,   notice that $\dim _\varphi$ depends on the ordering of
the isomorphism classes of  the simple right modules. Therefore when
we refer to a $\dim _\varphi$ function or to dimension vectors we
implicitly assume that we have chosen an ordering of the simple
modules. If we explicitly state that the semisimple artinian ring
$R$ is isomorphic to $M_{n_1}(D_1)\times \cdots \times
M_{n_k}(D_k)$, for $D_1,\dots ,D_k$ division rings, then we assume
we are choosing an ordered family of representatives of the
isomorphism classes of simple right (or left) $R$-modules
$(V_1,\dots,V_k)$ such that $\End _R(V_i)\cong D_i$ for $i=1,\dots
,k$.

To easy the work with the elements in $\No^*$ we shall use the
following definitions.

\begin{Def} \label{supports} Let $\mathbf{x}=(x_1,\dots ,x_k)\in (\No^*) ^k$. We define
\[\mathrm{supp}\, (\mathbf{x})=\{i\in \{1,\dots, k\}\mid x_i\neq 0\}\]
and we refer to this set as the \emph{support of $\mathbf{x}$}. We
also define \[\infsupp  (\mathbf{x})=\{i\in \{1,\dots, k\}\mid x_i=
\infty \},\]  we refer to this set as the \emph{infinite support of
$\mathbf{x}$}.
\end{Def}

\subsection{Passing modulo the Jacobson radical}

First we recall the following well known Lemma

\begin{Lemma} \label{full} Let $R$ be a ring with Jacobson radical $J(R)$. Let
$P$ and $Q$ be  projective right $R$-modules.
\begin{itemize}
\item[(i)] Assume $P$ and $Q$ are finitely generated. If there
exists a   projective right $R/J(R)$-module $X$ such that
$P/PJ(R)\cong Q/QJ(R)\oplus X$ then there exists a projective right
$R$-module $Q'$ such that $P\cong Q\oplus Q'$ and $Q'/Q'J(R)\cong
X$.

\item[(ii)] Assume   only that $Q$ is finitely generated. If $f\colon
P/PJ(R)\to Q/QJ(R)$ is an onto module homomorphism then $Q$ is
isomorphic to a direct summand of $P$.
\end{itemize}
\end{Lemma}

If, in the above Lemma,  neither $P$ nor $Q$ are finitely generated
then even the weaker divisibility property (ii) is lost. It was
shown in \cite{P1} that it is still true that projective
modules isomorphic modulo the Jacobson radical are isomorphic. We
recall this fundamental result in the next statement together with a
weaker property on lifting pure monomorphisms.

We recall that a right module monomorphism $f\colon M_1\to M_2$ is
said to be a pure monomorphism if, for any left module $N$,
$f\otimes _RN\colon M_1\otimes _R N\to M_2\otimes _RN$ remains a
monomorphism. For example, if $f$ is a (locally) split monomorphism
then it is pure. If $f$ is a monomorphism between two projective
modules then $f$ is pure if and only if $\mathrm{coker}\, f$ is a
flat module if and only if $f$ is locally split.

\begin{Th}\label{Pavel}  Let $R$ be any
ring, and let $P$ and $Q$ be projective right $R$-modules.
\begin{itemize}
\item[(i)] {\rm \cite[Proposition~6.1]{plans}} A module homomorphism $f\colon P\to Q$
is a pure monomorphism if and only if so is the induced map
$\overline{f}\colon P/PJ(R)\to Q/QJ(R)$.

\item[(ii)] {\rm \cite[Theorem~2.3]{P1}} If $f\colon P/PJ(R)\to Q/QJ(R)$ is an isomorphism of right
$R/J(R)$-modules then there exists an isomorphism of right
$R$-modules $g\colon P\to Q$ such that the induced morphism
$\overline g\colon P/PJ(R)\to Q/QJ(R)$ coincides with $f$.
\end{itemize}
\end{Th}

Theorem~\ref{Pavel}(ii) allows us to see the monoids of projective
modules over a ring $R$ as  submonoids of the monoids of projective
modules over $R/J(R)$. To give the assertion in a more precise way
we shall use the following notion (cf. \cite{FH2}).

\begin{Def} \label{deffull} A submonoid $A$ of a monoid $C$ is said to be a \emph{full
submonoid} of $C$ if for any $x\in A$ and any $t\in C$, $x+t\in A$
implies $t\in A$. If $f\colon A\to C$ is an injective monoid
homomorphism and $\im (f)$ is a full submonoid of $C$ we say that
$f$ is a \emph{full embedding}.

A \emph{full affine monoid} is a full submonoid of a finitely
generated free commutative monoid, and a \emph{full affine
embedding} is a full embedding into a finitely generated free
commutative monoid.
\end{Def}

See Proposition~\ref{relace} for a characterization of full affine
submonoids of $\No ^k$.

We note that in the terminology of \cite{bruns2} a full affine
embedding is a pure embedding of monoids.

\begin{Cor} \label{modjr} Let $R$ be a ring with Jacobson radical $J(R)$, and let $\pi \colon R\to R/J(R)$
denote the canonical projection.  Then:
\begin{itemize}
\item[(i)] $V(\pi)\colon V(R)\to V(R/J(R))$ is a full embedding of monoids with order
unit. In particular, the algebraic preorder on $V(\pi)(V(R))$
coincides with the one induced by the algebraic preorder on
$V(R/J(R))$.
\item[(ii)] $V^*(\pi)\colon V^* (R)\to V^* (R/J(R))$ is an
injective monoid morphism.
\end{itemize}
\end{Cor}

There is an interesting intermediate submonoid between $V(R)$ and
$V^*(R)$.

\begin{Def}\label{cuntzmonoid} Let $R$ be a ring. Set $W(R)=W(R_R)$
to be  the additive monoid of isomorphism classes of projective
right $R$-modules that are pure submodules of a finitely generated
free right $R$-module. The addition on $W(R)$ is induced by the
direct sum of modules.

Analogously, $W({}_RR)$ is the additive monoid of isomorphism
classes of projective left $R$-modules that are pure submodules of a
finitely generated free right $R$-module.
\end{Def}

For example, if $R=\mathcal{C}([0,1])$ is the ring of real valued
continuous functions defined on  the interval $[0,1]$ then the ideal
$$I=\{f\in R\mid \mbox{ there exists $\varepsilon >0$ such that }
f([0,\varepsilon])=0\}$$ is countably generated, projective and pure
inside $R$. Therefore $\langle I\rangle\in W(R)\setminus V(R)$.

The notation  $W(R)$ is borrowed from the $C^*$-algebra world, as we
think on $W(R)$ as the discrete analogue of the Cuntz monoid (cf.
\cite{elliott})

The following result, which is a consequence of
\cite[Theorem~7.1]{plans} and Theorem~\ref{Pavel}, describes one way
to obtain elements in $W(R)\setminus V(R)$ and which is the only one
when the ring $R$ is semilocal.

\begin{Prop} \label{wr} Fix $n\in \N$. Let $R$  be a  ring. Let $P_1$, $P_2$ be finitely generated projective right $R/J(R)$ modules
such that $(R/J(R))^n\cong P_1\oplus P_2$.  Then the following
statements are equivalent
\begin{itemize}
\item[(i)] There exists a projective right $R$-module $P$ such that
$P/PJ(R)\cong P_1$.
\item[(ii)] There exists a pure right submodule $M$ of $R^n$ such
that $M/MJ(R)\cong P_1$.
\item[(iii)] There exists a  projective left $R$-module
$Q$ such that $Q/J(R)Q\cong \mathrm{Hom}_{R/J(R)} (P_2, R/J(R))$.
\end{itemize}
In this case $P$ and $Q$ are countably generated pure submodules of
$R^n$. They are finitely generated if and only if there exists a
projective right $R$-module  $P'$ such that $P'/P'J(R)\cong P_2$.
\end{Prop}

Observe that, by Theorem~\ref{Pavel}, the isomorphism class of the
module $P$ appearing in Proposition~\ref{wr} is an element of
$W(R_R)$ and the isomorphism class of $Q$ gives an element of
$W({}_RR)$. Therefore if $P$ is not finitely generated, $\langle
P/PJ(R)\rangle\le \langle (R/J(R))^n\rangle$ in $W(R/J(R))$ but $P$
is not a direct summand of $R^n$, so that $\langle P\rangle $ and
$\langle R^n\rangle $ are uncomparable in $W(R)$. So that, in
general, the algebraic preorder on $W(R)$ does not coincide with the
order induced by the algebraic preorder on $W(R/J(R))$.

We will go back to the monoid $W(R)$ in \S \ref{final}, as we are
mainly interested in noetherian rings and then, clearly,
$V(R)=W(R)$. Results of Lazard \cite{lazard} show that this also
holds  just assuming ascending chain condition on annihilators.

\begin{Prop} \label{radicalfg} Let $R$ be a ring such that, for any $n\in \N$,
$M_n(R)$ has the ascending chain condition on right annihilators of
elements. Then a pure submodule of a finitely generated free right
$R$-module is finitely generated.
\end{Prop}

\begin{Proof} Combine \cite[Lemme 2(i)]{lazard} with the argument in
\cite[Corollary~3.6]{FHS}.
\end{Proof}

\subsection{Semilocal rings}

Let  $R$ be a  semilocal ring such that $R/J(R)\cong
M_{n_1}(D_1)\times \cdots \times M_{n_k}(D_k)$ for suitable
 division rings $D_1,\dots ,D_k$. Fix an onto ring homomorphism
$\varphi \colon R\to  M_{n_1}(D_1)\times \cdots \times M_{n_k}(D_k)$
such that $\mathrm{Ker}\, \varphi =J(R)$. Then there is an induced
ring isomorphism $\overline{\varphi}\colon R/J(R)\to
M_{n_1}(D_1)\times \cdots \times M_{n_k}(D_k)$, so that we have a
dimension function $\dim _{\overline{\varphi}}$, cf.
\S~\ref{semisimple}. For any countably generated projective right
$R$-module $P$, set
\[\dim _\varphi (\langle P\rangle): =\dim
_{\overline{\varphi}}(\langle P\otimes _RR/J(R)\rangle)=\dim
_{\overline{\varphi}}(\langle P/PJ(R)\rangle).\] By
Corollary~\ref{modjr}, $\dim _\varphi (V(R))$ is a full affine
submonoid of $\No ^k$ with order unit $(n_1,\dots ,n_k)$ and $\dim
_\varphi (V^*(R))$ is a submonoid of $(\No^*)^k$. For some
considerations on the order induced on $V^*(R)$ by $V^*(R/J(R))$ we
refer to Corollary~\ref{algorder} and to \S~\ref{final}.

It was shown in \cite{FH}
 that the full affine property characterizes the monoids $A$
with order unit that can be realized as  $V(R)$ of some semilocal
ring $R$. More precisely, if $A$ is a full affine submonoid of
$\No^k$ with order unit $(n_1,\dots ,n_k)$ then there exist a
semilocal hereditary ring $R$, $D_1,\dots ,D_k$ division rings and
an onto ring homomorphism $\varphi\colon R\to M_{n_1}(D_1)\times
\cdots \times M_{n_k}(D_k)$ with kernel $J(R)$ such that $\dim
_\varphi V(R)=A$.

Note that since  over a hereditary ring any projective module is a
direct sum of finitely generated projective modules, it follows that
for a hereditary ring $R$ as above
\[\dim _\varphi V^*(R)= A+\infty\cdot A\subseteq  (\N^*_0)^k\]
see Corollary~\ref{allfg} and Proposition~\ref{fas}.

\section{Semilocal rings: The noetherian case}\label{noetherian}

We start this section   recalling some results on projective modules
over noetherian semilocal rings from \cite{fair} and adapting them
to our purposes.  We  also state in \ref{main} our main
characterization theorem.

 It is well known that the trace ideal
  of a projective module is an idempotent ideal. Whitehead in
\cite{whitehead} characterized idempotent ideals that are trace
ideals of countably generated projective modules. His results yield
that in a noetherian ring any idempotent ideal is a trace ideal of a
countably generated projective module.  In \cite{fair}, P\v r\'\i
hoda noted that Whitehead's ideas can be extended to prove that if
$I$ is an idempotent ideal of a noetherian ring $R$ then any
finitely generated projective $R/I$-module can be extended to a
projective $R$-module. For further quotation we state these results.

\begin{Prop} \label{symtraces} Let $R$ be a  noetherian ring. Then the
following statements are equivalent for a two sided ideal $I$
\begin{itemize}
\item[(i)] $I^2=I$.
\item[(ii)] There exists a countably generated projective right
$R$-module $P$ such that $\mathrm{Tr}\,(P)=I$.
\item[(iii)] For any finitely generated projective right
$R/I$-module $P'$ there exists a countably generated projective
right $R$-module $P$ such that $P/PI\cong P'$ and $I\subseteq
\mathrm{Tr}\,(P)$.
\item[(iv)] There exists a countably generated projective left
$R$-module $Q$ such that $\mathrm{Tr}\,(Q)=I$.
\item[(v)] For any finitely generated projective left
$R/I$-module $Q'$ there exists a countably generated projective left
$R$-module $Q$ such that $Q/IQ\cong Q'$ and $I\subseteq
\mathrm{Tr}\,(Q)$.
\end{itemize}
\end{Prop}

\begin{Proof} Combine \cite[Corollary~2.7]{whitehead} with \cite[Lemma~2.6]{fair}.
\end{Proof}

Trace ideals of projective modules keep memory of the semisimple
factors of the projective module.

\begin{Lemma} \label{simplefactors} Let $R$ be a semilocal ring, and let $P$ be a   projective
right module with trace ideal $I$. Let $V_R$ be a simple right
$R$-module with endomorphism ring $D$, and let
$W=\mathrm{Hom}_D(V,D)\cong \mathrm{Hom}_R(V,R/J(R))$ be its dual
simple left $R$-module. Then the following statement are equivalent:
\begin{itemize}
\item[(i)] $V$ is a quotient of $P$.
\item[(ii)] $V$ is a quotient of $I$.
\item[(iii)] $I+r_R(V)=I+l_R(W)=R$.
\item[(iv)] $W$ is a quotient of $I$.
\end{itemize}
In particular, if $I$ is also the trace ideal of a left projective
module $Q$ then the above statements are also equivalent to the
fact that $W$ is a quotient of $Q$.
\end{Lemma}

\begin{Proof} The equivalence of (i) and (ii) is a particular case
of \cite[Lemma~3.3]{fair}. It is clear that (iii) is equivalent to
(ii) because, for a semilocal ring, $r_R(V)=l_R(W)$ is a maximal
two-sided ideal of $R$. Statements (iii) and (iv) are equivalent by
the symmetry of (iii).
\end{Proof}

\begin{Th}\label{behaviour} \emph{\cite{fair}} Let $R$ be a noetherian semilocal ring. Let $V_1,\dots,
V_k$ be an ordered set of representatives of  the isomorphism
classes of  simple right $R$-modules.

For $i=1,\dots ,k$, let $D_i=\mathrm{End}_R(V_i)$ and
$W_i=\mathrm{Hom}_{D_i}(V_i,D_i)\cong \mathrm{Hom}_{R}(V_i,R/J(R))$.
So that $W_1,\dots ,W_k$ is an ordered set of representatives of the
isomorphism classes of simple left $R$-modules. Let $S$ be a subset
of $\{1,\dots, k\}$.  Assume that there exists a countably generated
projective right $R$ module $P$ such that
\[P/PJ(R)\cong \left(\oplus _{i\in \{1,\dots ,k\}\setminus S}V_i^{n_i}\right) \oplus
\left(\oplus _{i\in S}V_i^{(\omega)} \right),\] where $n_i\in \No$.
Then the following statements hold:
\begin{itemize}
\item[(1)] There exists a countably generated projective right
$R$-module $P'$ such that $P'/P'J(R)\cong \oplus _{i\in
S}V_i^{(\omega)}$. Hence $P\cong P\oplus P'$.
\item[(2)] Let $I$ be the trace ideal of $P'$. Then $P/PI$ is a
finitely generated right $R/I$-module such that
\[P/PI\otimes _{R/I}\left(R/I\right)/J(R/I)\cong P/P(I+J(R))\cong \oplus _{i\in \{1,\dots ,k\}\setminus S}V_i^{n_i}.\]
\item[(3)] There exists a countably generated projective left
$R$-module $Q$ such that
\[Q/J(R)Q\cong \left(\oplus _{i\in \{1,\dots ,k\}\setminus S}W_i^{n_i}\right) \oplus \left(\oplus _{i\in
S}W_i^{(\omega)} \right)\]
\item[(4)] There exists a countably generated projective left
$R$-module $Q'$ such that $Q'/J(R)Q'\cong \oplus _{i\in
S}W_i^{(\omega)}$. Hence $Q\cong Q\oplus Q'$.
\end{itemize}
Therefore, $V^*(R_R)\cong V^*({}_RR)$ and,   fixing $\varphi \colon
R\to M_{n_1}(D_1)\times \cdots \times M_{n_k}(D_k)$ an onto ring
homomorphism with kernel $J(R)$, we obtain that $\dim _{\varphi}
V^*(R_R)= \dim _{\varphi}V^*({}_RR)$ and that $\dim _\varphi
V(R)=\left(\dim _{\varphi} V^*(R)\right)\cap \No^k$.
\end{Th}

\begin{Proof}$(1).$ The existence of $P'$ follows from
\cite[Proposition~3.4]{fair}. The isomorphism $P\cong P\oplus P'$
follows from Theorem~\ref{Pavel}(ii).

Statement $(2)$ is also part of \cite[Proposition~3.4]{fair}.

By Proposition~\ref{symtraces}, $I$ is also the trace ideal of a
projective left $R$-module $M$. As $M^{(\omega)}/J(R)M^{(\omega)}$
is semisimple and contains all semisimple factors of $M$ it follows
from Lemma~\ref{simplefactors} that
$M^{(\omega)}/J(R)M^{(\omega)}\cong \oplus _{i\in
S}(V^*_i)^{(\omega)}$. Therefore taking $Q'=M^{(\omega)}$ we deduce
that the first statement of $(4)$ holds.

By $(2)$, $P/PI$
 is a finitely generated $R/I$-module. Therefore $\overline Q=\mathrm{Hom}_{R/I}(P/PI,
 R/I)$ is a finitely generated projective left $R/I$-module such
 that $\overline{Q}/J(\overline{Q})\cong \oplus _{i\in
 \{1,\dots ,k\}\setminus S}W_i^{n_i}$. By Proposition~\ref{symtraces}, there exists a
 projective left $R$-module $Q_1$ such that $Q_1/IQ_1\cong
 \overline{Q}$. Then $Q=Q_1\oplus Q'$ fulfills the requirements of
 statement $(3)$ and the second half of statement $(4)$.

Finally, note that the assignment $\langle P\rangle \mapsto \langle
Q\rangle$ induces an isomorphism between $V^*(R_R)$ and $V^*({}_RR)$
such that $\dim _{\varphi}(V^*(R_R))=\dim _{\varphi}(V^*({}_RR))$.
The claim on $\dim_\varphi V(R)$ follows either from (2) or from
combining Proposition~\ref{wr} with Proposition~\ref{radicalfg}.
\end{Proof}

In Example~\ref{noniso} we will see that the isomorphism between
$V^*({}_RR)$ and $V^*(R_R)$ does not exist for general semilocal
rings.

As a corollary of Theorem~\ref{behaviour} we note that, in the
context of noetherian semilocal rings, the  divisibility property of
Lemma~\ref{full}(ii) still holds for general projective modules.

\begin{Cor}\label{weakdiv} Let $R$ be a noetherian semilocal ring. Let $P$ and $Q$ be projective right $R$-modules
such that   $P/PJ(R)$ is isomorphic to a direct summand of $Q/QJ(R)$
then $P$ is isomorphic to a direct summand of $Q$.
\end{Cor}

\begin{Proof} Since any projective module is a sum of countably
generated projective modules we may assume that $P$ and $Q$ are
countably generated \cite[Proposition~2.50]{libro}.

Let $V_1,\dots, V_k$ be an ordered set of representatives of  the
isomorphism classes of  simple right $R$-modules. Since $P/PJ(R)$ is
a homomorphic image of $Q/QJ(R)$ there exist $S'\subseteq S\subseteq
\{1,\dots ,k\}$ such that
\[Q/QJ(R)\cong \left(\oplus _{i\in \{1,\dots ,k\}\setminus S}V_i^{n_i}\right) \oplus
\left(\oplus _{i\in S}V_i^{(\omega)} \right)\]
and
\[ P/PJ(R)\cong
\left(\oplus _{i\in \{1,\dots ,k\}\setminus S'}V_i^{m_i}\right)
\oplus \left(\oplus _{i\in S'}V_i^{(\omega)}\right)
\]
where $n_i$ and $m_j$ are in $\No$, and $n_i-m_i\in \No$ for any
$i\in \{1,\dots ,k\}\setminus S$.

By Theorem~\ref{behaviour}, there exists a countably generated
projective module $Q'$ such that $Q'/Q'J(R)\cong \oplus _{i\in
S}V_i^{(\omega)}.$ Let $I$ be the trace ideal of $Q'$. Again by
Theorem~\ref{behaviour}, $\overline{P}=P/PI$ and $\overline{Q}=Q/QI$
are finitely generated projective right $R/I$-modules.

Now
\[\overline{Q}/\overline{Q}J(R/I)\cong Q/Q(I+J(R))\cong \oplus _{i\in \{1,\dots ,k\}\setminus
S}V_i^{n_i}\] and
\[\overline{P}/\overline{P}J(R/I)\cong P/P(I+J(R))\cong \oplus _{i\in \{1,\dots ,k\}\setminus
S}V_i^{m_i}\] By Corollary~\ref{modjr}, there exist a finitely
generated projective right $R/I$-module $\overline{X}$ such that
$\overline{X}/\overline{X}J(R/I)\cong \oplus _{i\in \{1,\dots
,k\}\setminus S}V_i^{n_i-m_i}.$ By Proposition~\ref{symtraces},
there exists a countably generated projective right $R$-module $X$
such that $X/XI\cong \overline{X}$. By Theorem~\ref{Pavel}(ii),
$Q\cong P\oplus Q'\oplus X$.
\end{Proof}

An alternative way to state   Corollary~\ref{weakdiv} is thinking on
$V^*(R)$ as  a  monoid with the partial order induced by the
algebraic order of $V^*(R/J(R))$.

\begin{Cor}\label{algorder} Let $R$ be a   semilocal ring.
Consider the following relation over $V^*(R)$, $\langle P\rangle \le
\langle Q\rangle $ if and only if $P/PJ(R)$ is isomorphic to a
direct summand of $Q/QJ(R)$. Then
\begin{itemize}
\item[(i)] $\langle P\rangle \le
\langle Q\rangle $ if and only if there exists a pure embedding
$f\colon P\to Q$.
\item[(ii)] $\le $ is a partial order relation that refines the algebraic order on $V^*(R)$.
\item[(iii)] If, in addition, $R$ is noetherian then the partial
order induced  by $\le $ over $V^*(R)$ is the algebraic order.
\end{itemize}
\end{Cor}

\begin{Proof} $(i).$ If  $\langle P\rangle \le
\langle Q\rangle $ then there exists a splitting monomorphism
$\overline{f}\colon P/PJ(R)\to Q/QJ(R)$ which by
Theorem~\ref{Pavel}(i) lifts to a pure monomorphism $f\colon P\to
Q$.  Conversely, if $f\colon P\to Q$ is a pure monomorphism of right
$R$-modules then the induced map $f\otimes _R R/J(R)\colon P\otimes
_R R/J(R)\to Q \otimes _R R/J(R)$ is a pure monomorphism of
$R/J(R)$-modules. Since $R/J(R)$ is semisimple, $f\otimes _RR/J(R)$
is a split monomorphism.

$(ii).$ It is clear that $\le $ is reflexive and transitive. As it
is already observed in $\cite{P1}$, Theorem~\ref{Pavel} implies that
$\le$  is also antisymmetric.

If $P$ is isomorphic to a direct summand of $Q$, then $P/PJ(R)$ is
also isomorphic to a direct summand of $Q/QJ(R)$. Hence $\langle
P\rangle \le \langle Q\rangle$, that is, $\le$ refines the algebraic
order on $V^*(R)$.

$(iii).$ It is a consequence of Corollary~\ref{weakdiv}.
\end{Proof}

After some amount of work, it will turn out that
Proposition~\ref{symtraces} and Theorem~\ref{behaviour} contain all
the information needed to describe $V^* (R)$ for $R$ a noetherian
semilocal ring.

\begin{Def} \label{defequations} Let $k\ge 1$. A submonoid $M$ of $(\No^*)^k$ is said to be a
monoid \emph{defined by a system of equations} provided that there
exist $D\in M_{n\times k}(\No)$, $E_1,$ $E_2\in M_{\ell \times
k}(\No)$ and $m_1,\dots ,m_n\in \N$ , $m_i\ge 2$ for any $i\in
\{1,\dots ,n\}$, such that $M$ is the set of   solutions in
$(\No^*)^k$ of the system of equations
\[D\left(\begin{array}{c}t_1\\\vdots \\ t_k\end{array}\right)\in \left(\begin{array}{c}m_1\No^*\\
\vdots \\ m_n\No^*\end{array} \right) \quad(*) \qquad \mbox{and} \qquad
E_1\left(\begin{array}{c}t_1\\\vdots \\
t_k\end{array}\right)=E_2\left(\begin{array}{c}t_1\\\vdots \\
t_k\end{array}\right) \quad(**)\] where $\ell$, $n\ge 0$. By convention,
$\ell$ or $n$ equal to zero means that either $(*)$ or $(**)$ are
empty systems.
\end{Def}

As we shall recall in Proposition~\ref{relace},   any full affine
monoid of $\No ^k$ is of the form $M\cap \No ^k$ where $M$ is a
submonoid of $(\No^*)^k$  defined by a system of equations.

Now we can state our main theorem,

\begin{Th} \label{main} Let $k\in \N$. Let $M$ be a submonoid of $(\No^*)^k$
containing  $(n_1,\dots ,n_k)\in \N ^k$. Then the following
statements are equivalent:
\begin{itemize}
\item[(1)] $M$ is defined by a system of equations.
\item[(2)] For any field $F$ there exist a noetherian semilocal
$F$-algebra $R$, a semisimple $F$-algebra $S=M_{n_1}(D_1)\times
\dots \times M_{n_k}(D_k)$, where $D_1,\dots,D_k$ are division
rings, and an onto morphism of $F$-algebras $\varphi \colon R\to S$
with $\mathrm{Ker}\, \varphi =J(R)$ such that $\dim _\varphi
 V^*(R) =M$. In particular, $\dim _\varphi  V(R) =M\cap \No^k$.
\item[(3)] There exist a noetherian semilocal
ring $R$, a semisimple ring $S=M_{n_1}(D_1)\times \dots \times
M_{n_k}(D_k)$, were $D_1,\dots,D_k$ are division rings, and an onto
ring morphism $\varphi \colon R\to S$ with $\mathrm{Ker}\, \varphi
=J(R)$ such that $\dim _\varphi  V^*(R) =M$. Therefore, $\dim
_\varphi  V(R) =M\cap \No^k$.
\end{itemize}
\end{Th}

\begin{Remark} We follow the notation of
Definition~\ref{defequations}. As it is done for full affine monoids
in \cite[Exercise~6.4.16]{BH} or \cite[Proof of
Theorem~2.29]{bruns2}, if $M\subseteq (\No^*)^k$ is defined by a
system of equations  as in \ref{defequations} then it is isomorphic
to the submonoid $M'$ of $(\No^*)^{k+n}$ defined by  system of
linear diophantine equalities
\[D\left(\begin{array}{c}t_1\\\vdots \\ t_k\end{array}\right)=\left(\begin{array}{ccc}m_1&\cdots &0\\
\vdots &\ddots &\vdots \\ 0&\cdots&m_n \end{array}\right) \left(\begin{array}{c}t_{k+1}\\
\vdots \\ t_{k+n}\end{array} \right) \qquad \mbox{and} \qquad
E_1\left(\begin{array}{c}t_1\\\vdots \\
t_k\end{array}\right)=E_2\left(\begin{array}{c}t_1\\\vdots \\
t_k\end{array}\right)\] The isomorphism is given by the assignment
\[(x_1,\dots ,x_k)\mapsto (x_1,\dots ,x_k, \frac 1{m_1}\sum
_{i=1}^kd_{1j}x_j,\dots ,\frac 1{m_n}\sum _{i=1}^kd_{nj}x_j),\]
where we make the convention $\frac \infty{m_i}=\infty$.

Therefore it is important to take into account that we are
considering our monoids always inside some fixed $(\No^*)^k$ or, in
the ring context, that we are viewing $V^*(R)$ as a submonoid of
$V^*(R/J(R))$.
\end{Remark}

The monoid $\No^*$ is not cancellative, therefore the solutions of
two systems of equations may coincide over $\No$ but be different
when considered over $\No ^*$. We illustrate this phenomena with an
easy example.

\begin{Ex} The set of solutions of the equation $x=y$ in $\No ^2$ is
$M=\{(n,n)\mid n\in \No\}$, and the set of solutions in $(\No^*)^2$
is $M+\infty\cdot M=M\cup \{(\infty,\infty) \}$.

The set of solutions of $2x=x+y$ in $\No ^2$ is, of course, also $M$
but in $(\No^*)^2$  is $M_1=M\cup \{(\infty, n)\mid n\in \No^*\}$.

Finally, the set of solutions of $2x+y=x+2y$ in $(\No^*)^2$  is
$M_1\cup \{(n,\infty)\mid n\in \No^*\}$.
\end{Ex}

 Theorem~\ref{main} shows that, for noetherian semilocal rings, the
 description of $V^*(R)$ viewed inside $V^*(R/J(R))$ nicely extends the one of $V(R)$ inside $V(R/J(R))$ (cf.
 Proposition~\ref{relace}). In \S \ref{final} we give examples showing
 that the picture for general semilocal rings must be more
 complicated.

 \section{Semilocal principal ideal domains} \label{PID}

We recall that a ring $R$ is a \emph{principal ideal domain} if $R$
is a right and left principal ideal domain, that is, if every right
ideal of $R$ has the form $aR$ for some $a\in R$ and every left
ideal of $R$ has the form $Ra$ for some $a\in R$.

Semilocal principal ideal domains are a source of semilocal
noetherian rings such that all projective modules are free. Our aim
in this section is to construct semilocal PID's with certain types
of semisimple factors.

Let $R$ be a commutative ring. Let $k\ge 1$, and let
$\mathcal{M}_1,\dots ,\mathcal{M}_k$ be  different maximal ideals of
$R$.  The
 localization of $R$ at the set $\Sigma =R\setminus
\left(\mathcal{M}_1\cup \dots \cup \mathcal{M}_k\right)$ is a
semilocal ring such that modulo its Jacobson radical is isomorphic
to $R/\mathcal{M}_1\times \cdots \times R/\mathcal{M}_k$.

Fuller and Shutters  observed in \cite{FS} that  the same procedure
to construct semilocal rings can be extended to, non necessarily
commutative, principal ideal domains by using Ore  localization.

\begin{Prop}\emph{\cite[Proposition~4]{FS}} \label{pid} Let $\varphi \colon R\to
S$ be a surjective ring homomorphism of a principal ideal domain $R$
onto a semisimple artinian ring $S$. Let $\Sigma =\{a\in R\mid
\varphi (a)\mbox{ is invertible in }S\}$. Then:
\begin{itemize}
\item[(i)] $\Sigma$ is a right and left Ore set.
\item[(ii)] The Ore localization $R_\Sigma$ of $R$ with
respect to $\Sigma$ is a semilocal principal ideal domain, and the
extension $\overline{\varphi}\colon R_\Sigma \to S$ of $\varphi$
induces an isomorphism $R_\Sigma /J(R_\Sigma)\cong S$.
\end{itemize}
\end{Prop}

Next result gives a source of examples where to apply
Proposition~\ref{pid}.

Let $E$ be any ring, and let $\alpha\colon E\to E$ a (unital) ring
morphism. The skew polynomial ring or the twisted polynomial ring is
the ring
\[E[x;\alpha]=\{p(x)=p_0x^m+\cdots +p_m\mid m\in \No \mbox{ and }p_i\in E\mbox{ for }
i=0,\dots,m\}\] with componentwise addition and multiplication
induced by the rule $xr=\alpha{(r)}x$ for any $r\in E$.

It is well known that if $E$ is a division ring and $\alpha$ is an
automorphism then $E[x;\alpha]$ has a right and a left division
algorithm, hence, it is a principal ideal domain.

\begin{Prop} \label{skew} Let $E$ be a field. Let $\alpha\colon E\to E$ be a
field automorphism of order $n$ with fixed field $E^\alpha=\{a\in E
\mid \alpha(a)=a\}$. Then the skew polynomial ring $R=E[x;\alpha ]$
has a simple factor isomorphic to $M_n(E^\alpha)$.

Moreover, if $E$ is infinite then, for any $k\in \N$, $R$ has a
factor isomorphic to $M_n(E^\alpha)^k$.
\end{Prop}

\begin{Proof} We may assume that $n>1$.

Note the following  fact that will be useful throughout the proof:

(*)\emph{ Let $p(x)=x^m+p_1x^{m-1}+\cdots +p_m\in R$ be such that
$p_m\neq 0$. If $a\in E$ satisfies that $ap(x)\in p(x)R$ then
$\alpha ^m(a)=a$.}

As $\alpha$ has order $n$, the center of $R$ contains (in fact
coincides with) $E^\alpha[x^n]$. Therefore, for any $0\neq t\in
E^\alpha$, the right ideal $(x^n-t)R$ is  two-sided. As $R$ is a
right principal ideal domain, (*) yields that $(x^n-t)R$ is a
maximal two-sided ideal, so that $R/(x^n-t)R$ is  a simple artinian
ring. We claim that if  $t=r^n$ for some $r\in E^\alpha$ then
$R/(x^n-t)R\cong M_n(E^\alpha)$. To prove this we need to find a
simple right $R/(x^n-t)R$-module such that its endomorphism ring is
$E^\alpha$ and its dimension over $E^\alpha$ is $n$.

In $E^\alpha[x]\subseteq R$ we have a decomposition $x^n-r^n=(x-r
)q(x)$. As $x^n-r^n$ is central in $R$, $V=R/(x-r)R$ is a right
$R/(x^n-r^n)R$-module.  It is readily checked that $V$ is a right
$E$-vector space of dimension $1$, therefore it is a simple right
$R/(x^n-r^n)R$-module.

$\mathrm{End}_R(V)=\mathcal{I}/(x-r)R$, where $\mathcal{I}=\{p(x)\in
R\mid p(x)(x-r)\in (x-r)R\}$ is the idealizer of $(x-r)R$ in $R$. As
any  $p(x)\in R$ can be written in a unique way as $a+(x-r)q(x)$ for
$a\in E$,
\[\mathrm{End}_R(V)\cong \mathcal{I}\cap E=\{a\in E\mid a(x-r)\in (x-r)R\}\]
By (*), $ \mathcal{I}\cap E=E^\alpha$ and, hence,
$\mathrm{End}_R(V)\cong E^\alpha$. Since ${}_{E^\alpha} V\cong
{}_{E^\alpha} E$ and, by Artin's Theorem, $[E:E^\alpha]=n$ we deduce
that the dimension of $V$ over its endomorphism ring is $n$ as
desired.

Now assume that $E$, and hence $E^\alpha$, is infinite. Fix $k\in
\N$. Let $r_1,\dots ,r_k\in E^\alpha$ be such that $r_1^n,\dots
,r_k^n$ are $k$ different elements. Consider the ring homomorphism
\[\Phi \colon R\to R/(x^n-r_1^n)R\times \cdots \times
R/(x^n-r_k^n)R\] defined by $\Phi(p(x))=(p(x)+(x^n-r_1^n)R, \dots
,p(x)+(x^n-r_k^n)R)$. Clearly, $\mathrm{Ker}\, \Phi =\cap
_{i=1}^k(x^n-r_i^n)R$. Since, for any $i=1,\dots,k$,
$(x^n-r_i^n)R+\bigcap _{j\neq i}(x^n-r_j^n)R=R$ we deduce that
$\Phi$ is also onto. Therefore, by the first part of the proof,
$\Phi$ induces an isomorphism $R/\mathrm{Ker}\, \Phi\cong
M_n(E^\alpha)^k$.
\end{Proof}

In Theorem~\ref{pbrealization} we will use the following examples.

\begin{Exs}\label{mnk} Let $n,k\in \N$.
\begin{itemize}
\item[(i)] There exists a semilocal $\Q$-algebra $R$ that is a principal
ideal domain such that $R/J(R)\cong M_n(\Q)^k$.
\item[(ii)] Let $F$ be any field, and consider the transcendental
extension of $F$, $E=F(t_1,\dots ,t_n)$. Let $\alpha\colon E\to E$
be the automorphism of $E$ that fixes $F$ and satisfies that $\alpha
(t_i)=t_{i+1}$ for $i=1,\dots ,n-1$ and $\alpha (t_n)=t_1$. Then
there exists a semilocal $F$-algebra $R$, that is a principal ideal
domain, such that $R/J(R)\cong M_n(E ^\alpha)^k$.
\end{itemize}
In both cases all projective right or left modules over $R$ are
free.
\end{Exs}

\begin{Proof} (i) Let $\Q \subseteq E$ be a Galois field extension
with Galois group $G\cong \Z/n\Z$. Let $\alpha \colon E\to E$ be a
generator of $G$. By Proposition~\ref{skew}, there exists an onto
ring homomorphism $\varphi \colon E[x;\alpha]\to M_n(\Q)^k$. By
Proposition~\ref{pid},   $\Sigma =\{a\in R\mid \varphi (a)\mbox{ is
invertible in }M_n(\Q)^k\}$ is a right and left Ore set and
$(E[x;\alpha])_\Sigma$  has the desired properties.

(ii) Proceed as in (i) combining Proposition~\ref{skew} with
Proposition~\ref{pid}.
\end{Proof}

\section{Pullbacks of rings} \label{pullbacks}

We shall use ring pullbacks to construct noetherian semilocal rings
with prescribed $V^*(R)$. In this section we study when ring
pullbacks are semilocal and noetherian. We start fixing some
notation.

\begin{Not} \label{pb} Let $R_1$, $R_2$ and $S$ be rings with ring homomorphisms
$j_i\colon R_i \to S$, for $i=1,2$. Let $R$ be the pullback of these
rings. That is, $R$ fits into the pullback diagram
$$\begin{CD}
R_1@> j_1 >>S\\
@A i_1 AA @AA j_2 A\\
R@>> i_2> R_2
\end{CD}$$
So that it can be described as $$R=\{(r_1,r_2)\in R_1\times R_2\mid
j_1(r_1)=j_2(r_2)\}$$ and the maps $i_1$ and $i_2$ are just the
canonical projections. \end{Not}

A ring homomorphism $\varphi \colon R\to S$ is said to be
\emph{local} if, for any $r\in R$, $\varphi (r)$ is a unit of $S$ if
and only if $r$ is a unit of $R$. The following  deep result by Rosa
Camps and Warren Dicks (see \cite[Theorem~1]{campsdicks} or
\cite[Theorem~4.2]{libro}) characterizes semilocal rings in terms of
local morphisms.

\begin{Th}\label{local} A ring $R$ is semilocal if and only if it has a local ring homomorphism into a
semilocal ring.\end{Th}

We note the following elegant  corollary of Theorem~\ref{local}.

\begin{Cor}\label{localpb} In the situation of Notation~\ref{pb}, $R$ is a local
subring of $R_1\times R_2$. In particular, the pullback  of two
semilocal rings is a  semilocal ring.
\end{Cor}

\begin{Proof} Note that if $(r_1,r_2)\in R\subseteq R_1\times R_2$ is a unit of $R_1\times R_2$ then its inverse
$(r_1^{-1},r_2^{-1})$ also satisfies that
$j_1(r_1^{-1})=j_2(r_2^{-1})$. Hence $(r_1^{-1},r_2^{-1})\in R$, and
we deduce that the inclusion $R\to R_1\times R_2$ is a local ring
homomorphism.
\end{Proof}

\begin{Lemma}\label{subring} Let $T$ be a subring of a ring $R$, and
assume that there exists a two-sided ideal $I$ of $R$ such that
$I\subseteq T$ and that $R/I$ is finitely generated as a left
$T/I$-module. Then:
\begin{itemize}
\item[(i)] ${}_TR$ is finitely generated.
\item[(ii)] If $R$ and $T/I$ are left noetherian rings then so is $T$.
\end{itemize}
\end{Lemma}

\begin{Proof} $(i)$ Let $x_1,\dots ,x_n \in R$ be such that
$x_1+I,\dots ,x_n+I$ generate $R/I$ as a left $T/I$-module. Then
$R=Tx_1+\cdots +Tx_n+T\cdot 1$ is finitely generated as a left
$T$-module.

$(ii)$ Let $J$ be a left ideal of $T$. As $IJ$ is a left ideal of
$R$, it is finitely generated as a left $R$-module. By (i), it is
finitely generated as a left $T$-module.

The left $R$-module $RJ/IJ$ is also a finitely generated left
$R/I$-module, hence it is a noetherian left $T/I$-module. Therefore
$J/IJ\subseteq RJ/IJ$ is a finitely generated left $T/I$-module.
Since
\[0\to IJ\to J\to J/IJ\to 0\]
we can conclude that ${}_TJ$ is finitely generated. This proves that
$T$ is left noetherian.
\end{Proof}

\begin{Prop}\label{noetherianpb} In the situation of Notation~\ref{pb}, assume that
 $j_1$  is surjective and that ${}_{R_2}S$ is finitely generated. If $R_i$ is a left noetherian
ring, for $i=1,2$, then $R$ is left noetherian.
\end{Prop}

\begin{Proof} As
$j_1$ is onto, $i_2$ is also an onto ring homomorphism with kernel
$\mathrm{Ker}\, j_1\times \{0\}$. Let $T=i_1(R)$, and note that
$\mathrm{Ker}\, j_1$ is a two-sided ideal of $R_1$ that is contained
in $T$. As a first step we shall prove that $T$ is left noetherian
and that ${}_TR_1$ is finitely generated. In view of
Lemma~\ref{subring} we only need to prove that ${}_TS\cong
R_1/\mathrm{Ker}\, j_1$ is finitely generated.

By assumption, there exist $s_1,\dots ,s_n\in S$ such that $S=\sum
_{i=1}^nR_2s_i$. Fix an element $s\in S$, there exist $r_2^1,\dots
r_2^n\in R_2$ such that $s=\sum _{i=1}^nr_2^i\cdot s_i=\sum
_{i=1}^nj_2(r_2^i)s_i$. Since $j_1$ is onto, for $i=1,\dots ,n$,
there exists $r_1^i\in R_1$ such that $j_1(r_1^i)=j_2(r_2^i)$. Hence
$s=\sum _{i=1}^nj_1(r_1^i)s_i=\sum _{i=1}^nr_1^i\cdot s_i$. Since
$(r_1^i,r_2^i)\in R$, $r_1^i\in T$ for $i=1,\dots ,n$. This shows
that $S=\sum _{i=1}^nTs_i$, so that ${}_TS$ is finitely generated.

We want to prove that any left ideal of $R$ is finitely generated.
Let $I$ be a left ideal of $R$ contained in $\mathrm{Ker}\,
j_1\times \{0\}$. Hence $I=I_1\times 0$ with $I_1$ a left ideal of
$T$, as $T$ is left noetherian $I$ is finitely generated.

Now, let $I$ be any left ideal of $R$. Since $i_2$ is onto and $R_2$
is left noetherian, $i_2(I)$ is a left ideal of $R_2$, finitely
generated by elements $r_2^1,\dots ,r_2^n\in R_2$ say. Fix
$r_1^1,\dots ,r_1^n\in R_1$ such that $(r_1^i,r_2^i)\in I$. If $x\in
I$ then there exists $y_1,\dots ,y_n\in R$ such that $x-\sum
_{i=1}^ny_i(r_1^i,r_2^i) \in I\bigcap\left(\mathrm{Ker}\, j_1\times
\{0\}\right)$. Therefore $I=I\bigcap\left(\mathrm{Ker}\, j_1\times
\{0\}\right)+\sum _{i=1}^nR(r_1^i,r_2^i)$. By the previous case,
$I\bigcap\left(\mathrm{Ker}\, j_1\times \{0\}\right)$ is finitely
generated, therefore $I$ is finitely generated.
\end{Proof}

In the next result we compute the Jacobson radical for some
pullbacks of rings.

\begin{Lemma} \label{radicalpb} In the situation of Notation~\ref{pb}, assume that
 $j_1$  is an    onto ring homomorphism  such that
 $\mathrm{Ker}\, j_1 \subseteq J(R_1)$ and $j_1(J(R_1))\supseteq
 j_2(J(R_2))$. Then $J(R)$ fits into the induced pullback diagram
 $$\begin{CD}
J(R_1)@> j_1 >>j_1(J(R_1))\\
@A i_1 AA @AA j_2 A\\
J(R)@>> i_2> J(R_2)
\end{CD}$$
and $R/J(R)\cong R_2/J(R_2)$. In particular, if
$J(R_1)=\mathrm{Ker}\, j_1$ and $J(R_2)=\mathrm{Ker}\, j_2$ then
$J(R)=J(R_1)\times J(R_2)$.
\end{Lemma}

\begin{Proof} Let $J$ be the pullback of the induced maps $j_1\colon
J(R_1)\to j_1(J(R_1))$ and $j_2\colon J(R_2)\to j_1(J(R_1))$. Since,
by Corollary~\ref{localpb}, $R$ is a local subring of $R_1\times
R_2$ it follows that $J\subseteq J(R)$.

In order to prove the reverse inclusion consider $(r_1,r_2)\in
J(R)$. Being $j_1$ onto, $i_2$ is also onto, hence $r_2\in J(R_2)$.
Since $j_1(r_1)=j_2(r_2)\subseteq j_1(J(R_1))$, we deduce that
$r_1\in \mathrm{Ker}\, j_1+J(R_1)=J(R_1)$. Therefore $(r_1,r_2)\in
J$.

Since $\mathrm{Ker}\, i_2=\mathrm{Ker}\, j_1\times \{0\}\subseteq
J(R)$ and $i_2$ is onto, it immediately follows that the map
\[R\stackrel{i_2}\to R_2\stackrel{\pi}\to R_2/J(R_2),\]
where $\pi$ denotes the canonical projection, induces an isomorphism
between $R/J(R)$ and $R_2/J(R_2)$.

The claim when $J(R_i)=\mathrm{Ker}\, j_i$ follows from the fact
that the pullback of   zero homomorphisms is the product.
\end{Proof}

Milnor proved in \cite{milnor}  that, in the situation of
Notation~\ref{pb} and assuming that $j_1$ is onto, the category of
projective modules over $R$ is the \emph{pullback} of the category
of projective modules over $R_1$ and of the category of projective
modules over $R_2$. We are interested in the following version of
this result,

\begin{Th}\label{milnor} \emph{\cite[Theorem~2.1]{milnor}} In the situation of Notation~\ref{pb}, assume that
 $j_1$  is an    onto ring homomorphism.  Then there are  pullbacks of
 monoid morphisms
 $$\begin{CD}
V(R_1)@> V(j_1) >>V(S)\\
@A V(i_1) AA @AA V(j_2) A\\
V(R)@>> V(i_2)> V(R_2)
\end{CD}\qquad\mbox{ and }\qquad \begin{CD}
V^*(R_1)@> V^*(j_1) >>V^*(S)\\
@A V^*(i_1) AA @AA V^*(j_2) A\\
V^*(R)@>> V^*(i_2)> V^*(R_2)
\end{CD}$$
\end{Th}

\begin{Rem}{\rm
In fact, Milnor proved that having projective modules $P_1$ over
$R_1$ and $P_2$ over $R_2$ such that $P_1 \otimes_{R_1} S$ and $P_2
\otimes_{R_2}S$ are isomorphic $S$-modules, then there exists a
(unique) projective $R$-module $P$ such that $P_1 \simeq P
\otimes_{R} R_1$ and $P_2 \simeq P \otimes_{R} R_2$. Moreover, if
$P_1$ is finitely generated over $R_1$ and $P_2$ is finitely
generated over $R_2$, then $P$ is finitely generated over $R$. He
did not state explicitly that $P$ is countably generated $R$-module
if $P_1$ is countably generated over $R_1$ and $P_2$ is countably
generated $R_2$. However, it follows directly from his proof.}
\end{Rem}

Now we state the precise result we will be using in \S
\ref{construction}.

\begin{Cor} \label{milnorsemilocal} In the situation of Notation~\ref{pb}, let $R_1$
and $R_2$ be  semilocal rings, and let  $S$ be a semisimple artinian
ring. Fix an isomorphism $\varphi\colon S\to M_{m_1}(D'_1)\times
\cdots \times M_{m_\ell}(D'_\ell)$ for suitable division rings
$D'_1,\dots ,D'_\ell$. Assume that $j_1$ is an onto ring
homomorphism with kernel $J(R_1)$, and that $J(R_2)\subseteq
\mathrm{Ker} \, j_2$.  If $R_2/J(R_2)\cong M_{n_1}(D_1)\times \cdots
\times M_{n_k}(D_k)$ for $D_1,\dots ,D_k$ division rings, and $\pi
\colon R_2\to M_{n_1}(D_1)\times \cdots \times M_{n_k}(D_k)$ is an
onto morphism with kernel $J(R_2)$ then
\begin{itemize}
\item[(i)] $i_2$ induces an onto ring homomorphism
$\overline{i_2}\colon R\to M_{n_1}(D_1)\times \cdots \times
M_{n_k}(D_k)$ with kernel $J(R)$. In particular, $R$ is a semilocal
ring and  $R/J(R)\cong R_2/J(R_2)$.
\item[(ii)] Let $\alpha \colon \dim _{\overline{i_2}} \colon V^*(R_2)\to (\No^*)^{\ell}$ be the monoid homomorphism
induced by $i_2$ and $\varphi$. Then $$\dim _{\overline{i_2}}
V^*(R)=\{x\in \dim _{\pi }V^*(R_2)\subseteq (\No^*)^k\mid \alpha
(x)\in \dim _{\varphi j_1}V^*(R_1)\}.$$
\end{itemize}
Moreover, if $R_1$ and $R_2$ are noetherian, and $S$  is finitely
generated, both as a left and as a right $j_2(R_2)$-module, then $R$
is noetherian and $\dim _{\overline{i_2}} (V(R))= \dim
_{\overline{i_2}} V^*(R)\cap (\No )^k$.
\end{Cor}

\begin{Proof}  Statement $(i)$ follows from Lemma~\ref{radicalpb} and
Corollary~\ref{localpb}. Statement $(ii)$ is a consequence of
Theorem~\ref{milnor}.

The final part of the Corollary follows from
Proposition~\ref{noetherianpb} and the fact that over a noetherian
semilocal ring a projective module is finitely generated if and only
if it is finitely generated modulo the Jacobson radical (cf.
Proposition~\ref{radicalfg} or Theorem~\ref{behaviour}).
\end{Proof}

We single out the following particular case of
Corollary~\ref{milnorsemilocal}.

\begin{Cor}\label{milnorintersection} In the situation of Notation~\ref{pb}, for $i=1,2$,
assume that $j_i$ is onto and $\mathrm{Ker}\, j_i=J(R_i)$. Assume
$S=M_{n_1}(D_1)\times \cdots \times M_{n_k}(D_k)$ where $D_1,\dots
,D_k$ are division rings. Then
\begin{itemize}
\item[(i)] $i_2$ (and $i_1$) induces an onto ring homomorphism
$\overline{i}\colon R\to M_{n_1}(D_1)\times \cdots \times
M_{n_k}(D_k)$ with kernel $J(R)$.
\item[(ii)]  $\dim _{\overline{i}} V^*(R)=\dim _{j_1}V^*(R_1)\cap \dim
_{j_2}V^*(R_2)$.
\end{itemize}
\end{Cor}

\section{Noetherian semilocal rings with prescribed $V^*(R)$}
\label{construction}

Now we have all the elements to construct noetherian semilocal rings
with prescribed $V^*(R)$ and to prove the realization part of
Theorem~\ref{main}. We explain the basic constructions in the
following two examples.

\begin{Ex}\label{congruences}  Let $k, m\in \N$, and let
$a_1,\dots,a_k\in \No$. Assume $(n_1,\dots ,n_k)\in \N ^k$ is such
that $a_1n_1+\cdots +a_kn_k= m \ell\in \N$. Let $F$ be a field, and
let $F\subseteq F_2$ be a field extension such that there exist a
semilocal principal ideal domain $R_1$, that is also an $F$-algebra,
with $R_1/J(R_1)\cong M_m(F_2)$. Then for any intermediate field
$F\subseteq F_1\subseteq F_2$ such that $[F_2:F_1]<\infty$ there
exist a noetherian semilocal $F$-algebra $R$ and an onto ring
homomorphism $\varphi\colon R\to M_{n_1}(F_1)\times \cdots \times
M_{n_k}(F_1)$ with $\mathrm{Ker}\, \varphi =J(R)$ such that $\dim
_\varphi V^*(R)$ is exactly the set of solution in $(\No^*)^k$ of
the congruence $a_1t_1+\cdots +a_kt_k\in  m\No^*$.
\end{Ex}

\begin{Proof}
 Fix
$j_1\colon M_\ell (R_1)\to M_{m\ell}(F_2)$ an onto ring homomorphism
with kernel $J(M_\ell (R_1))=M_\ell (J(R_1))$.

Set $R_2=M_{n_1}(F_1)\times \cdots \times M_{n_k}(F_1)$, and
consider the morphism
\[\begin{array}{cccc}j_2\colon &R_2&\longrightarrow
&M_{m\ell}(F_2)\\ &(r_1,\dots ,r_k)&\mapsto &
\left(\begin{array}{ccc}\smallmatrix r_1&\cdots &0\\
\vdots &\ddots\! ^{a_1)} &\vdots \\ 0&\cdots& r_1\endsmallmatrix\
&\cdots &
0\\
&\ddots & \\
0&\cdots &\smallmatrix r_k&\cdots &0\\
\vdots&\ddots\! ^{a_k)} &\vdots \\ 0&\cdots& r_k\endsmallmatrix
\end{array}\right)
\end{array}\]
Note that $(V(R_2),\langle R_2\rangle) \cong (\No^k,(n_1,\dots
,n_k))$ and $V^*(R_2)\cong (\No^*)^k$; $V(M_{m\ell}(F_2))$ is
isomorphic to the monoid $\No$ with order unit $m\cdot \ell$ and
$V^*(M_{m\ell}(F_2))\cong (\No)^*$. Then $j_2$ induces the morphism
of monoids $f \colon (\No^*)^k\to \No^*$ defined by $f(x_1,\dots
,x_k)=a_1x_1+\cdots +a_kx_k$, cf. \S \ref{semisimple}.

 Let $R$ be the ring
defined by the pullback diagram
$$\begin{CD}
M_\ell(R_1)@> j_1 >>M_{m\ell}(F_2)\\
@A i_1 AA @AA j_2 A\\
R@>> \varphi > M_{n_1}(F_1)\times \cdots \times M_{n_k}(F_1)
\end{CD}$$
 Since $[F_2:F_1]<\infty$ we can apply  Corollary~\ref{milnorsemilocal} to deduce $R$ is a noetherian semilocal
 $F$-algebra
 and that $\varphi$ is an onto ring homomorphism with kernel $J(R)$.

 Now we
 compute $\dim _{\varphi}V^*(R)$ using also Corollary~\ref{milnorsemilocal}. We have chosen $R_1$  such that $\dim
_{j_1}V^*(M_\ell(R_1))=m\No^*$.  Therefore $(x_1,\dots ,x_k)\in \dim
_{\varphi}V^*(R)$ if and only if $f(x_1,\dots, x_k)=a_1x_1+\dots
+a_kx_k\in m\No^*$. That is, $\dim _{\varphi}V^*(R)$ is exactly the
set of solutions in $(\No^*)^k$
 of the congruence $a_1t_1+\cdots+a_kt_k\in m\No^*$ as desired.
\end{Proof}

\begin{Ex}\label{equalities}  Let $k\in \N$, and let
$a_1,\dots,a_k, b_1,\dots ,b_k\in \No$. Let  $(n_1,\dots ,n_k)\in \N
^k$ be such that $a_1n_1+\cdots +a_kn_k= b_1n_1+\cdots +b_kn_k\in
\N$. For any field extension $F\subseteq F_1$, there exist a
noetherian semilocal $F$-algebra $R$ and an onto ring homomorphism
$\varphi \colon R\to M_{n_1}(F_1)\times \cdots \times M_{n_k}(F_1)$
with kernel $J(R)$ such that $\dim _{\varphi}V^*(R)$ is the set of
solutions in $(\No^*)^k$ of the equation $a_1t_1+\cdots +a_kt_k=
b_1t_1+\cdots +b_kt_k$.
\end{Ex}

\begin{Proof}Set $m=a_1n_1+\cdots +a_kn_k= b_1n_1+\cdots +b_kn_k$.

Let $R_1$ be a noetherian semilocal $F$-algebra such that
$R_1/J(R_1)\cong F_1\times F_1$, and all projective modules over
$R_1$ are free. For example, we could take the commutative principal
ideal domain $R_1\cong F_1[x]_\Sigma$ with $\Sigma=(F_1[x])\setminus
\left(xF_1[x]\cup (x-1)F_1[x]\right).$

Let $j_1\colon M_m(R_1)\to M_m(F_1)\times M_m(F_1)$ be an onto ring
homomorphism with kernel $J(M_m(R_1))$. Set $R_2=M_{n_1}(F_1)\times
\dots \times M_{n_k}(F_1)$. Consider the morphism $j_2\colon
R_2\longrightarrow M_{m}(F_1)\times M_m(F_1)$ defined by
\[j_2(r_1,\dots ,r_k)=
\left( \left(\begin{array}{ccc}\smallmatrix r_1&\cdots &0\\
\vdots &\ddots\! ^{a_1)} &\vdots \\ 0&\cdots& r_1\endsmallmatrix
&\cdots &
0\\
&\ddots & \\
0&\cdots &\smallmatrix r_k&\cdots &0\\
\vdots&\ddots\! ^{a_k)} &\vdots \\ 0&\cdots& r_k\endsmallmatrix
\end{array}\right),\left(\begin{array}{ccc}\smallmatrix r_1&\cdots &0\\
\vdots &\ddots\! ^{b_1)} &\vdots \\ 0&\cdots& r_1\endsmallmatrix
&\cdots &
0\\
&\ddots & \\
0&\cdots &\smallmatrix r_k&\cdots &0\\
\vdots&\ddots\! ^{b_k)} &\vdots \\ 0&\cdots& r_k\endsmallmatrix
\end{array}\right)\right)
\]
Note that $j_2$ induces the morphism of monoids $f\colon
(\No^*)^k\to \No^*\times \No^*$ defined by $f(x_1,\dots
,x_k)=(a_1x_1+\cdots +a_kx_k, b_1x_1+\cdots +b_kx_k)$, cf. \S
\ref{semisimple}. Hence, $f(n_1,\dots ,n_k)=(m,m)$.

Let $R$ be the ring defined by the pullback diagram
$$\begin{CD}
M_m(R_1)@> j_1 >>M_{m}(F_1)\times M_m(F_1)\\
@A i_1 AA @AA j_2 A\\
R@>> \varphi > M_{n_1}(F_1)\times \cdots \times M_{n_k}(F_1)
\end{CD}$$

 By
Corollary~\ref{milnorsemilocal}, $R$ is a noetherian semilocal
$F$-algebra
 and $\varphi$ is an onto ring homomorphism with kernel $J(R)$.
We have chosen $R_1$  such that $\dim
_{j_1}V^*(M_\ell(R_1))=\{(x,x)\mid x\in \No^*\}$. Also by
Corollary~\ref{milnorsemilocal},
 $(x_1,\dots ,x_k)\in \dim _{\varphi}V^*(R)$ if and only
if $f(x_1,\dots ,x_k)\in \dim _{j_1}V^*(M_\ell(R_1))$ if and only if
$a_1x_1+\dots +a_kx_k=b_1x_1+\cdots+b_kx_k$ as desired.
\end{Proof}

\begin{Th}\label{pbrealization} Let $k\ge 1$, and let $F$ be a field. Let $M$ be a submonoid of  $(\N _0^*)^k$
defined by a system of equations and containing an element
$(n_1,\dots ,n_k)\in \N ^k$. Then there exist a noetherian semilocal
$F$-algebra $R$, a semisimple $F$-algebra $S=M_{n_1}(E)\times \dots
\times M_{n_k}(E)$, where $E$ is a field extension of $F$,  and an
onto morphism of $F$-algebras $\varphi \colon R\to S$ with
$\mathrm{Ker}\, \varphi =J(R)$ satisfying that $\dim _\varphi
V^*(R)=M$.
\end{Th}

\begin{Proof} Let  $M$ be defined by the
system of equations,
\[D\left(\begin{array}{c}t_1\\\vdots \\ t_k\end{array}\right)\in \left(\begin{array}{c}m_1\No^*\\
\vdots \\ m_n\No^*\end{array} \right)\quad(*)\qquad \mbox{and}
\qquad
E_1\left(\begin{array}{c}t_1\\\vdots \\
t_k\end{array}\right)=E_2\left(\begin{array}{c}t_1\\\vdots \\
t_k\end{array}\right)\quad(**)\] where  $D\in M_{n\times k}(\No)$,
$E_1,$ $E_2\in M_{\ell \times k}(\No)$, $n,\ell\ge 0$ and $m_1,\dots
,m_n\in \N$ , $m_i\ge 2$ for any $i\in \{1,\dots ,n\}$.

\noindent\textbf{Step 1.} \emph{There exist a field $E$ containing
$F$, a noetherian semilocal $F$-algebra $R_1$ and an onto ring
homomorphism $\varphi _1\colon R_1\to M_{n_1}(E)\times \cdots \times
M_{n_k}(E)$ such that $\dim _{\varphi _1} V^*(R_1)$ is the set of
solution in $(\No^*)^k$ of the system of congruences $(*)$.}

If $n=0$, that is, if $(*)$ is empty we set $E=F$,
$R_1=M_{n_1}(E)\times \cdots \times M_{n_k}(E)$ and $\varphi
_1=\mathrm{Id}$. Assume $n>0$, therefore we may also assume that all
the rows of $D$ are non zero.

Consider the field extension $F\subseteq F'=F(t^i_j\mid i=1,\dots
,n,\, j=1,\dots ,m_i)$. For each $i=1,\dots ,n$ consider the
automorphism $\alpha _i$ of $F'$ that fixes $F_i=F(t^s_j\mid s\neq
i)\subseteq F'$,  maps $t^i_j$ to $t^i_{j+1}$ for $1\le j<m_i$, and
maps $t^i_{m_i}$ to $t^i_1$. Note that $\alpha _i$ has order $m_i$.
Let $G$ be the group of permutations of $m_1+\cdots +m_k$ variables.
Then $G$ acts on $F'$. Set $E=(F')^G$, and note that $E\subseteq F'$
is a finite field extension.

By Example~\ref{mnk}(ii), for each $i=1,\dots ,n$, we can construct
a principal ideal domain  such that modulo the Jacobson radical is
isomorphic to $M_{m_i}((F')^{\alpha _i})$.  By
Example~\ref{congruences}, for $i=1,\dots, n$, there exist a
noetherian semilocal $F$-algebra $L_i$  and an onto ring morphism
$\pi _i\colon L_i \to M_{n_1}(E)\times \dots \times M_{n_k}(E)$ with
kernel $J(L_i)$ and such that $\dim _{\pi _i}V^*(L_i)$ is the set of
solutions in $(\No^*)^k$ of the $i$-th congruence in $(*)$.

Let $R_1$ be the pullback of $\pi _i$, $i=1,\dots ,n$. By
Corollary~\ref{milnorsemilocal},
  $R_1$ is a noetherian
semilocal $F$-algebra. By Corollary~\ref{milnorintersection}, there
exist an onto ring homomorphism $\varphi_1\colon R_1\to
M_{n_1}(E)\times \dots \times M_{n_k}(E)$, with kernel $J(R_1)$,
such that $\dim _{\varphi _1} V^*(R_1)$ is the monoid  of solutions
of (*). This concludes the proof of the first step.

\medskip

\noindent\textbf{Step 2.} \emph{There exist a  noetherian semilocal
$F$-algebra $R_2$ and an onto ring homomorphism $\varphi _2\colon
R_2\to M_{n_1}(E)\times \cdots \times M_{n_k}(E)$ such that $\dim
_{\varphi _2} V^*(R_2)$ is the set of solution in $(\No^*)^k$ of the
system of equations $(**)$.}

If $\ell=0$, that is, if $(**)$ is empty we set
$R_2=M_{n_1}(E)\times \cdots \times M_{n_k}(E)$ and $\varphi
_2=\mathrm{Id}$. Assume $\ell >0$. Therefore,  we can assume that
none of the rows in $E_1$ and, hence, in $E_2$ are zero.

By Example~\ref{equalities}, for $i=1,\dots ,\ell$, there exist a
noetherian semilocal $F$-algebra $T_i$ and an onto ring morphism
$\pi _i\colon T_i \to M_{n_1}(E)\times \dots \times M_{n_k}(E)$ with
kernel $J(T_i)$ and such that $\dim _{\pi _i}V^*(T_i)$ is the set of
solutions in $(\No^*)^k$ of the $i$-th equation defined by the
matrices $E_1$ and $E_2$.

Let $R_2$ be the pullback of $\pi _i$, $i=1,\dots ,\ell$. By
Corollary~\ref{milnorsemilocal}, $R_2$ is a noetherian semilocal
$F$-algebra with an onto ring morphism $\varphi_2\colon R_2\to
M_{n_1}(E)\times \dots \times M_{n_k}(E)$ with kernel $J(R_2)$. By
Corollary~\ref{milnorintersection}, $\dim _{\varphi _2} V^*(R_2)$ is
the set of  solutions of $(**)$. This concludes the proof of step 2.

\medskip

Finally, set $R$ to be the pullback of $\varphi _i \colon R_i\to
M_{n_1}(E)\times \dots \times M_{n_k}(E)$, $i=1,2$. By
Corollary~\ref{milnorsemilocal}, $R$ is a noetherian semilocal
$F$-algebra with an onto ring morphism $\varphi\colon R\to
M_{n_1}(E)\times \dots \times M_{n_k}(E)$ with kernel $J(R)$. By
Corollary~\ref{milnorintersection}, the elements in $\dim _{\varphi
} V^*(R)$ are the solutions of $(*)$ and $(**)$.
\end{Proof}

\section{Solving systems in $\No$ and in $\No^*$: Supports of
solutions} \label{solving}

In this section we study the supports of the elements in monoids
defined by systems of equations. Our main aim is to show in
Proposition \ref{fas} that if $A$ is a full affine submonoid of
$\No^k$ then $A+\infty A$ is a monoid defined by a system of
equations.

Next result is quite easy but it is very important to keep it in
mind, for example in Definition~\ref{defsystems}. It shows that full
affine submonoids are closed by projections over the complementary
of supports of elements.

\begin{Prop}
\label{faprojections} Let $k\ge 1$. Let $A \subseteq \No ^k$ be a
full affine submonoid. Let $I \varsubsetneq \{1,\dots,k\}$ be the
support of an element of $A$, and denote by $p_I\colon \No ^k\to
\No^{\{1,\dots ,k\}\setminus I}$ the canonical projection. Then
$p_I(A)$ is a full affine submonoid of $\No^{\{1,\dots ,k\}\setminus
I}$.
\end{Prop}

\begin{Proof}
Clearly $p_I(A)$ is a  submonoid of $\No^{\{1,\dots ,k\}\setminus
I}$. We need to check the full affine property.

Let $a$, $b\in A$ be such that there exists $z\in \No^{\{1,\dots
,k\}\setminus I}$ satisfying that $p_I(a)+z=p_I(b)$.  Let $d\in A$
be such that $\supp (d)=I$. There exists $n\in \No$ such that
$a+c=b+nd$ for some $c\in \No ^k$. As $A$ is full affine, $c\in A$
and, hence, $z=p_I(c)\in p_I(A)$.
\end{Proof}

Let $C\subseteq \Q ^k$, and let $C^\perp=\{v\in \Q ^k\mid \langle
v,c\rangle=0\mbox{ for any }c\in C\}$ where $\langle -,-\rangle$
denotes the standard scalar product. If $X \subseteq \Q^k$, the
support of $X$ is defined by
$$\supp(X)=\bigcup _{x\in X}\supp (x)$$

 Let $A$ be a  submonoid of
$\No ^k$, and let  $B$ be the subgroup of $\Z ^k$ generated by $A$.
Then $A$ is full affine if and only if $B\cap \No ^k=A$ (see, for
example, \cite[Lemma~3.1]{FH}). So assume $A$ is full affine and
consider $B'=(A^\perp)^\perp$ which is a subgroup of $\Z^k$,
containing $B$,  defined by a system of diophantine linear equations
\[E\left(\begin{array}{c} t_1\\ \vdots \\ t_k \end{array}\right)=\left(\begin{array}{c}0\\ \vdots \\ 0
\end{array}\right)\]
with $E\in M_{(k-\ell) \times k}(\Z)$ and $\ell $ is the rank of the
group $B'$.

Since the rank of $B$ is also $\ell $, there exist $d_1,\dots,d_\ell
\ge 1$ such that $B'/B\cong \Z/d_1\Z\times \cdots \times \Z/d_\ell
\Z$. Equivalently, there exists a basis $v_1,\dots ,v_\ell$ of $B'$
such that $d_1v_1,\dots ,d_\ell v_\ell $ is a basis of $B$.
Therefore, an element   $x\in B'$
\[x=(x_1,\dots ,x_k)=\alpha _1v_1+\cdots +\alpha _\ell v_\ell \]
is in $B$ if and only if, for any $i=1,\dots ,\ell$, $\alpha _i\in
d_i\Z$. Since each $\alpha _i$ can be written as a  $\Q$-linear
combination of $x_1,\dots ,x_k$, by clearing denominators and
eliminating trivial congruences, we deduce that, for any  $x\in B'$,
 $x\in B$ if and only if it is a solution of
\[D\left(\begin{array}{c} t_1\\ \vdots \\ t_k \end{array}\right)\in \left(\begin{array}{c} m_1\Z\\ \vdots \\ m_n\Z
\end{array}\right)\]
where $0\le n \le \ell$, $D\in M_{n\times k}(\Z)$ and $m_i>1$ for
$i=1,\dots ,n$.

Adding to $D$ a suitable integral matrix of the form
$\left(\begin{smallmatrix}m_1&  \cdots &0\\ \vdots &\ddots  &\vdots
\\ 0& \cdots &m_n\end{smallmatrix}\right)\cdot A$, we can   also assume that $D\in M_{n\times
k}(\No)$.

In the next Proposition we collect the consequences of this
discussion.

\begin{Prop}\label{relace} \emph{(\cite[Exercise~6.4.16]{BH} or \cite[Proof of Theorem~2.29]{bruns2})} Let $k\ge 1$. Let $A$ be a full affine submonoid of
$\No ^k$, and let $\ell$ be the rank of the group generated by  $A$.
Then there exist $0\le n \le \ell$, $D\in M_{n\times k}(\No)$,
$E_1$, $E_2\in M_{(k-\ell) \times k}(\No)$ and $m_1,\dots ,m_n$
integers strictly bigger than one such that
\begin{itemize}
\item[(i)] $x=(t_1,\dots,t_k)\in \No^k$ is an element of $A$ if and only if it is a solution
of
\[D\cdot T\in \left(\begin{array}{c} m_1\No\\ \vdots \\ m_n\No
\end{array}\right)\qquad  \mbox{ and } \qquad E_1\cdot T=E_2\cdot T \]
where $T= (t_1,\dots ,t_k)^t$.

For $j=1,2$, let $r_i^j$ denote the $i$-th row of $E_j$. Then $E_1$
and $E_2$ can be chosen such that, for $i=1,\dots ,k-\ell$, $\supp
(r^1_i)\cap \supp (r^2_i)=\emptyset$.
\item[(ii)] The set of solutions of $E_1\cdot T=E_2\cdot T$ is
$A'=(A^\perp)^\perp \cap \No^k$.
\item[(iii)] There exists $d\in \N$ such that
$dA'\subseteq A$. In particular, \[\{I\subseteq \{1,\dots ,k\}\mid
\mbox{ there exists $a\in A$ such that
$\supp(a)=I$}\}=\]\[=\{I\subseteq \{1,\dots ,k\}\mid \mbox{ there
exists $a\in A'$ such that $\supp(a)=I$}\}.\]
\end{itemize}
\end{Prop}

\begin{Proof} Following with the notation in the remarks before Proposition~\ref{relace},
we can write the matrix $E=E_1-E_2$ where $E_1$ and $E_2$ are in
$M_{(k-\ell) \times k}(\No)$. Clearly, $E_1$ and $E_2$ can be
chosen in a way such that the $i$-th row of $E_1$ has disjoint
support with the $i$-th row of $E_2$. Then $(i)$ follows from the
fact that $A=B\cap \No ^k$.

The rest of the statement is clear.
\end{Proof}

Now we prove an auxiliary (and probably known) result that will be
useful to determine the supports of positive solutions of linear
diophantine equations.

\begin{Lemma} \label{lem}  Let $k\ge 1$, and  let $V$ be a subspace of
$\Q^k$. Then the following statements are equivalent,
\begin{itemize}
\item[(i)] $V^\perp \cap \N ^k\neq \emptyset$;
\item[(ii)] $\supp (V^\perp \cap \No ^k)=\{1,\dots
,k\}$;
\item[(iii)] $V\cap \N_0^k=\{0\}$.
\end{itemize}
\end{Lemma}

\begin{Proof} It is clear that (i) and (ii) are equivalent
statements and also that (i) implies (iii). We will show that (iii)
implies (ii).

The assumption in (iii) is equivalent to say that any element  $0\ne
v\in V$ has a component strictly bigger that zero an another one
strictly smaller that zero. As a first step we  show that $V^\perp$
cannot satisfy this condition.

Assume, by the way of contradiction, that $k$ is the minimal
dimension in which the conclusion of (ii) fails. So that there
exists $V\le \Q^k$  such that $V\cap \No^k= V^\perp \cap \No ^k=
\{0\}$. Note that $k$ and the dimension of $V$ must be strictly
bigger than $1$.

Let $v_1,\dots ,v_n$ be a basis of $V$ such that there exists $i\in
\supp (v_1)\setminus  \supp (\{v_2,\dots ,v_n\})$. Let $\pi \colon
\Q ^k \to \Q^{\{1,\dots, k\}\setminus \{i\}}$ denote the canonical
projection.

If $\pi (V)\cap \No^{\{1,\dots, k\}\setminus \{i\}}=\{0\}$ then, by
the minimality of $k$, there exists $v\in \Q ^k$ such that $0\ne
\pi(v)\in \pi (V)^\perp\cap \No^{\{1,\dots, k\}\setminus \{i\}}$.
Since $v$ can be chosen satisfying that $i\not \in \supp (v)$, we
would get $0\neq v\in V^\perp \cap  \N_0^k$, a contradiction. Let
$0\neq \lambda _1\pi (v_1)+\cdots +\lambda _n\pi (v_n)\in \pi
(V)\cap \No^{\{1,\dots, k\}\setminus \{i\}}$. Then $w=\lambda _1
v_1+\cdots +\lambda _nv_n\in V$. Since $V\cap \No^k=\{0\}$, $\lambda
_1\neq 0$.  Therefore, replacing $v_1$ by $w$ if necessary, we may
assume that $0\neq \pi(v_1)\in \No ^{\{1,\dots,k\}\setminus \{i\}}$ and that the
$i$-th component of $v_1$ is $<0$. Let $-a$ be such component.

Let $W$ be the subspace of $\Q ^{k\setminus\{i\}}$ generated by $\pi
(v_2),\dots ,\pi (v_n)$. Our hypothesis imply that $W\cap   \No
^{\{1,\dots ,k\}\setminus \{i\}}=\{0\}$. By the minimality of $k$,
there is $v\in \Q^k$ such that $0\neq \pi (v)\in W^\perp \cap \No
^{\{1,\dots ,k\}\setminus \{i\}}$. Therefore, as $b=\langle \pi
(v),\pi (v_1)\rangle \ge 0$, picking $v$ such that its $i$-th
component is $b/a$ we find that $0\neq av\in V^\perp \cap \No ^k$, a
contradiction. Therefore $V^\perp \cap \No ^k \neq \{0\}$ for any
$V$ such that $V\cap \No ^k =\{0\}$, as claimed.

Now assume that $V$ is a $\Q$-vector space satisfying (iii). Observe
first that $\supp (V^\perp)=\{1,\dots,k\}$, since otherwise there
would exist $i\in \{1,\dots,k\}$ such that the $i$-th component of
any element in $V^\perp$ is zero. Therefore, $e_i=(0,\dots,0,
1^{i)},0,\dots ,0)\in (V^\perp)^\perp=V$, which is a contradiction
with the assumption.

Let $I=\supp (V^\perp \cap \No ^k)$ and let $J=\{1,\dots
,k\}\setminus I$. We already know that $I\neq \emptyset$ and we want
to show $J=\emptyset$. Assume, by the way of contradiction, that
$J\neq \emptyset$. Set $\pi _I\colon \Q ^k\to \Q ^I$ and $\pi
_J\colon \Q ^k\to \Q ^J$ to be the canonical projections.

Pick $x\in V^\perp \cap \No ^k$ such that $\supp (x)=I$. Then for
any $v\in V^\perp$ there exists $n\in \No$ such that $\pi_I(nx+v)\ge
0$. By the definition of $I$, this implies that $\pi _J(v)=\pi
_J(nx+v)$ is either zero or it has a component $>0$ and another one
$<0$. Therefore, $\pi _J(V^\perp)\cap \No ^J=\{0\}$. By the first
part of the proof, there exists $w\in \Q ^k$ such that $0\neq \pi
_J(w)\in \pi _J(V^\perp)^\perp\cap \No ^J$.  Choosing $w$ such that
$\pi _I(w)=0$, we obtain that $0\neq w\in (V^\perp)^\perp\cap \No^k=
V\cap \No ^k$ which contradicts (iii). Therefore $J=\emptyset$.
\end{Proof}

Lemma~\ref{lem} yields a first characterization of the supports of
the elements in a full affine monoid. In following statements
if $X \subseteq (\No^*)^k$, $\Supp (X) \subseteq 2^{\{1,\dots,k\}}$ is the set
of supports of elements in $X$.

\begin{Cor} \label{lemsupports} Let $k\ge 1$.
Let $A $  be a full affine submonoid of $\No ^k$. Let $\emptyset
\neq I \subseteq \{1,\dots,k\}$, and denote by $\pi _I\colon \Q ^k
\to \Q ^I$  the canonical projection.
 Then  $I\in \Supp (A)$ if and only if $\pi _I
(A)^\perp\cap \No ^I=\{0\}$.
\end{Cor}

\begin{Proof} If $a\in A$ is such that $\supp (a)=I$ then, as  $\langle x,
\pi _I(a)\rangle=0$ for any $x\in \pi_I (A)^\perp$, it follows that
$\pi_I (A)^\perp\cap \No ^I=\{0\}$.

Conversely, if $\pi_I (A)^\perp\cap \No ^I=\{0\}$ then, by
Lemma~\ref{lem}, there exists $u\in (\pi_I (A)^\perp)^\perp\cap \N
^I$. Let $x=(x_1,\dots ,x_k)\in \No ^k$ be such that $\pi_I (x)=u$
and $x_i=0$ for any $i\in \{1,\dots ,k\}\setminus I$. Then $x\in
(A^\perp)^\perp\cap \No^k$. By Proposition \ref{relace} (iii), there
exists $d\in \N$ such that $a=dx\in A$. By construction, $\supp
(a)=I$.
\end{Proof}

A further characterization is the following.

\begin{Cor} \label{minimal} Let $k\ge 1$.
Let $A $  be a full affine submonoid of $\No ^k$. For any $\emptyset
\neq I \subseteq \{1,\dots,k\}$ denote by $\pi _I\colon \Q ^k \to \Q
^I$  the canonical projection. Let $v_1,\dots ,v_r\in \Q ^k$ be a
finite subset of $A^\perp$. Then, there exist $v_{r+1},\dots ,v_s\in
A^\perp$ such that $v_1,\dots v_s$ generate the $\Q$-vector space
$A^\perp$ and the set
\[\Scal (v_1,\dots ,v_s)=\{\emptyset \neq I\subseteq \{1,\dots
,k\}\mid\mbox{ for any $i=1,\dots ,s$, $\pi _I(v_i)$ is either
zero}\]\[\mbox{or it has a component $<0$ and a  component $>0$}\}\]
coincides with $\Supp (A\setminus \{0\})$.
\end{Cor}

\begin{Proof} Pick $v_{r+1},\dots ,v_s\in A^\perp$ such that  $v_1,\dots ,v_s$ generate
$A^\perp$ and  the set $\Scal (v_1,\dots ,v_s)$ has minimal
cardinality.

We claim that $v_1,\dots ,v_s$ have the desired property.  If
$\emptyset \neq I\subseteq \{1,\dots ,k\}$ is such that there exist
$a\in A$  such that $\supp (a)=I$ then $I\in \Scal (v_1,\dots ,v_s)$
because $\langle a,v_i\rangle =0$ for any $i=1,\dots ,s$.

Let $I\in \Scal (v_1,\dots ,v_s)$ and assume, by the way of
contradiction, that  $I\not \in \Supp (A)$. Then $\pi
_I(A^\perp)\cap \No ^I\neq \{0\}$, because otherwise
\[\pi _I(A)^\perp\cap \No ^I\subseteq \pi _I(A^\perp)\cap \No
^I=\{0\}\] and, by Corollary~\ref{lemsupports}, this implies that
there exists $a\in A$ with $\supp (a)=I$.

Let $v\in A^\perp$ be such that $0\neq \pi _I(v)\in\pi
_I(A^\perp)\cap \No ^I$. Notice that, $\Scal (v_1,\dots
,v_s,v)\subseteq \Scal (v_1,\dots ,v_s)$ and $I\in \Scal (v_1,\dots
,v_s)\setminus \Scal (v_1,\dots ,v_s, v)$. This contradicts the
minimality of the cardinality of $\Scal (v_1,\dots ,v_s)$. This
finishes the proof of the claim and of the Corollary.
\end{Proof}

Now we move to consider the solutions over   $(\No^*)^k$. In the
next lemma we see how to determine the set of   nonempty supports of
such monoids, which coincide with the set of infinite supports, for
the special kind of systems that appear in
Proposition~\ref{relace}(i).

\begin{Lemma} \label{infsupp} Let $k\ge 1$. Let $M$ be a submonoid of $(\No^*)^k$
defined by a system of equations
\[D\cdot T\in \left(\begin{array}{c}m_1\No^*\\
\vdots \\ m_n\No^*\end{array} \right)\quad (*)\qquad \mbox{and}
\qquad E_1\cdot T=E_2\cdot T\quad (**)\] where $D\in M_{n\times
k}(\No)$, $E_1,$ $E_2\in M_{\ell \times k}(\No)$ and $m_1,\dots
,m_n\in \N$ , $m_i\ge 2$ for any $i\in \{1,\dots ,n\}$.  For $j=1,2$
and $i=1,\dots ,\ell$, let $r_i^j$ denote the $i$-th row of $E_j$.
For $\emptyset \neq I\subseteq \{1,\dots ,k\}$, let $\pi _I\colon
(\No^* )^k\to (\No^*)^I$ denote the canonical projection. Then:
\begin{itemize}
\item[(i)] Let $N$ be the submonoid of $(\No^*)^k$ whose elements are the solutions
of the system of congruences $(*)$. Then, for any $i\in  \{1,\dots
,k\}$, the element $(0,\dots,0 ,\infty,0,\dots ,0)\in N$.
\item[(ii)] If $x\in M$ then also $\infty\cdot x\in M$.
\item[(iii)] If $x\in M$ then the element $x^*\in (\No^*)^k$ uniquely determined
by the property $\supp x^*=\infsupp x^*=\infsupp x$ also belongs to
$M$.
\item[(iv)]  Assume
that, for $i=1,\dots ,\ell$, $\supp(r^1_i)\cap
\supp(r^2_i)=\emptyset$, and let $\emptyset \neq I \subseteq
\{1,\dots ,k\}$ then there exists $x\in M$ such that $\supp (x)=I$
if and only if,
 for any $i=1,\dots
,\ell$,  $\pi _I(r^1_i-r^2_2)$ is either $0$ or it has a component
$>0$ and another one $<0$.
\end{itemize}
\end{Lemma}

\begin{Proof}
Statement (i) is trivial, and it allows us to prove the rest of the
statement just for $M$ defined by the system of linear diophantine
equations $(**)$.

Let  $x\in M$.  Fix $i\in \{1,\dots ,\ell\}$,  there are three
possible situations. The first one is  $0=\langle x,r_i^1\rangle
=\langle x,r_i^2\rangle$ which happens if and only if, for $j=1,2$,
$\supp(x)\cap \supp (r^j_i)=\emptyset$. The second one is  $0\neq
\langle x,r_i^1\rangle =\langle x,r_i^2\rangle\in \No$ which happens
if and only if, for $j=1,2$, $\infsupp (x)\cap \supp
(r_i^j)=\emptyset$ but $\supp (x)\cap \supp (r_i^j)\neq \emptyset$.
Finally,  $\langle x,r_i^1\rangle = \langle x, r_i^2 \rangle=\infty$
if and only if $\infsupp(x)\cap \supp(r_i^j)\neq \emptyset$ for
$j=1,2$. Then, in the three situations, it also follows that
$\langle \infty\cdot x,r_i^1\rangle =\langle \infty\cdot
x,r_i^2\rangle$ and $\langle x^*,r_i^1\rangle =\langle
x^*,r_i^2\rangle$. This shows that (ii) and (iii) hold.

To prove statement (iv) assume that $\supp (r^1_i)\cap \supp
(r^2_i)=\emptyset$ for $i=1,\dots ,\ell$. Let $\emptyset\neq
I\subseteq \{1,\dots ,k\}$ have the property required in the
statement. Let $x\in (\No^*)^k$ be such that $\supp (x)=\infsupp
(x)=I$. If $i\in \{1,\dots ,\ell\}$ is such that  $\pi
_I(r^1_i-r^2_i)$ is zero then $0=\langle x,r_i^1\rangle =\langle
x,r_i^2\rangle$, if $\pi _I(r^1_i-r^2_i)$ has a positive component
and a negative component then $\infty=\langle x,r_i^1\rangle
=\langle x,r_i^2\rangle$. This shows that $x$ satisfies $(**)$,
therefore it is an element of $M$.

To prove the converse, let $x\in M$. By $(ii)$ we may assume that
$x=\infty\cdot x$. Let $I=\supp (x)$, then one can proceed as in the
proof of $(ii)$ and $(iii)$ to show that $I$ has the required
property.
\end{Proof}

\begin{Prop}\label{fas} Let $k\ge 1$, and let $A$ be a full affine submonoid of $\No ^k$.
Then $$M=A + \{\infty\cdot a\mid a\in A\}$$ is a submonoid of ${(\No
^*)}^k$ defined by a system of equations.
\end{Prop}

\begin{Proof} We divide the proof into a couple of steps.

\noindent\textbf{Step 1.} \emph{Let $k\ge 1$, and let $A$ be a full
affine submonoid of $\No ^k$. Then there exists a submonoid  $M'$ of
${(\No ^*)}^k$ defined by a system of equations such that $M'\cap
\No^k=A$ and if $x\in M'$ then there exists $a\in A$ such that
$\infty\cdot x=\infty\cdot a\in M'$.}

By Proposition~\ref{relace} there exist $D\in M_{n\times k}(\No)$,
$E_1,$ $E_2\in M_{\ell \times k}(\No)$ and $m_1,\dots ,m_n\in \N$,
$m_i\ge 2$ for any $i\in \{1,\dots ,n\}$ such that the elements of
$A$ are the solutions of the system
\[D\left(\begin{array}{c}t_1\\\vdots \\ t_k\end{array}\right)\in \left(\begin{array}{c}m_1\No^*\\
\vdots \\ m_n\No^*\end{array} \right)\quad (*)\qquad \mbox{and}
\qquad
E_1\left(\begin{array}{c}t_1\\\vdots \\
t_k\end{array}\right)=E_2\left(\begin{array}{c}t_1\\\vdots \\
t_k\end{array}\right)\quad (**).\] For $j=1,2$, let $r_i^j$ denote
the $i$-th row of $E_j$. By Proposition~\ref{relace}, we can assume
that, for $i=1,\dots ,\ell$, $\supp (r_i^1)\cap \supp
(r_i^1)=\emptyset$. Set $v_1=r^1_1-r^2_1,\dots ,v_\ell =r^1_\ell
-r^2_\ell$. Notice that $v_1,\dots ,v_\ell\in A^\perp$ and, in fact,
generate $A^\perp$.

By Corollary~\ref{minimal}, there exist $v_{\ell+1},\dots ,v_s\in
A^\perp$ such that the set of supports of nonzero elements in $A$
coincides  with
\[\Scal (v_1,\dots ,v_s)=\{\emptyset \neq I\subseteq \{1,\dots
,k\}\mid\mbox{ for any $i=1,\dots ,s$, $\pi _I(v_i)$ is either
zero}\]\[\mbox{or it has a component $<0$ and a  component $>0$}\}\]
Now, for $i=\ell+1,\dots ,s$, write $v_i=r^1_i-r^2_i$ where $r_i^j
\in \No^k$ and $\supp (r_i^1)\cap \supp (r_i^1)=\emptyset$. For,
$j=1,2$, let $F_j$ be the matrix whose $i$-th row is $r^j_{\ell
+i}$.  Now, add to the initial system defining $A$ the equations
defined by
\[F_1\left(\begin{array}{c}t_1\\\vdots \\
t_k\end{array}\right)=F_2\left(\begin{array}{c}t_1\\\vdots \\
t_k\end{array}\right).\]

Let $M'$ be the set of solutions in $(\No^*)^k$, and note that $A$
is still the set of solutions in $\No ^k$ of this system so that
$M'\cap \No^k=A$. By Lemma~\ref{infsupp}(iv), the set of supports of
elements of $M'$ is exactly $\Scal (v_1,\dots ,v_s)$ which, by
construction, coincides with the set of supports of elements in $A$.
This implies that if $x\in M'$ then there exists $a\in A$ such that
$\infty\cdot x=\infty\cdot a$. By Lemma~\ref{infsupp}, if $x\in M'$
then $\infty\cdot x\in M'$. This finishes the proof of the first
step.

\noindent\textbf{Step 2.} \emph{The monoid $M$ in the statement is
defined by a system of equations}

Let \[\Scal =\{I\varsubsetneq \{1,\dots ,k\}\mid \mbox{there exists
$a\in A$ such that $\supp\, (a)=I$}\}.\] Notice that $\emptyset \in
\Scal$.

For any $I \in \Scal$, denote by $p_I\colon \No ^k\to \No^{\{1,\dots
,k\}\setminus I}$ the canonical projection. If $I\in \Scal$ then, by
Proposition~\ref{faprojections}, $A_I=p_I(A)$ is a full affine
submonoid of $\No^{\{1,\dots ,k\}\setminus I}$.

By step 1, there is a monoid $M'_I\subseteq (\No^*)^{\{1,\dots
,k\}\setminus I}$ defined by a system of equations and such that
$\No^{\{1,\dots ,k\}\setminus I}\cap M'_I=A_I$ and if $x\in M'_I$
then there exists $a\in A_I$ such that $\infty\cdot x=\infty\cdot
a$.

Set $M_I=\pi _I^{-1}(M'_I)$. Notice that $M_I$ is defined by the
same   system of equations defining $M'_I$ but considered over
$\No^k$.  Notice that $x \in \No^k\cap M_I$ if and only if
$\pi_I(x)\in A_I$. In particular, $A=\No^k\cap M_\emptyset$ and
$A\subseteq M_I$ for any $I\in \Scal$.

Since for any $I\in \Scal$, $M_I$ is defined by a system of
equations so is $\cap _{I\in \Scal}M_I$. We claim that $M=\cap
_{I\in \Scal}M_I$. We already know that
 $A\subseteq \cap _{I\in \Scal}M_I$, so that, by
 Lemma~\ref{infsupp}(iv), $M\subseteq \cap _{I\in \Scal}M_I$. To
 prove the other inclusion, let $x\in \cap _{I\in \Scal}M_I$. Let
 $I_1=\infsupp\, (x)$ and consider the element $x^*\in (\No^*)^k$ such that
 $I_1=\infsupp\, (x^*)=\supp (x^*)$. By
 Lemma~\ref{infsupp}(iii), $x^*\in \cap
_{I\in \Scal}M_I$. By construction, there exists $a\in A$ such that
$x^*=\infty\cdot a$, therefore $a\in M$. Since $\supp (a)=I_1$, we
deduce that $I_1\in \Scal$. Therefore, $\pi _{I_1}(x)=\pi
_{I_1}(a_1)$ for $a_1\in A$. This implies that $x=x^*+a_1$, so that
$x\in M$. This finishes the proof of the claim and the proof of the
Proposition.
\end{Proof}

\section{Systems of supports} \label{systems}
In order to conclude the proof of Theorem~\ref{main},  we need to
show that the monoids that appear as $V^*(R)$ for  noetherian
semilocal rings are defined by a system of equations. To this aim we
abstract the following class of submonoids of ${(\No ^*)}^k$.

\begin{Def} \label{defsystems}
Fix $k \in \mathbb{N}$ and an order unit $(n_1,\dots,n_k) \in
\mathbb{N}^k$. A \emph{system of supports}
$\mathcal{S}(n_1,\dots,n_k)$ consists of a collection $\cal S$ of
subsets of $\{1,\dots,k\}$ together with a system of monoids $A_I,
I \in \cal S$ such that the following conditions hold
\begin{enumerate}
\item[(i)] $\emptyset \in \cal S$, and $(n_1,\dots,n_k) \in A_{\emptyset}$.
\item[(ii)] For any $I \in \cal S$ the monoid $A_I$ is a  submonoid of
$\mathbb{N}_0^{\{1,\dots,k\} \setminus I}$. The monoid
$A_{\{1,\dots,k\}}$ is considered as a trivial monoid.
\item[(iii)]  $\mathcal{S}$ is closed under unions, and if $x\in A_I$ for some $I\in \mathcal{S}$ then $I\cup \supp(x)\in \cal S$.
\item[(iv)] Suppose that $I,K \in \cal S$ are  such that $I \subseteq K$ and
let $p \colon \mathbb{N}_0^{\{1,\dots,k\}\setminus I} \to
\mathbb{N}_0^{\{1,\dots,k\} \setminus K}$ be the canonical
projection. Then $p(A_I) \subseteq A_K$.
\end{enumerate}

If, for any $I\in \Scal$, the submonoids $A_I$ are full affine
submonoids of $\No ^{\{1,\dots,k\} \setminus I}$ then
$\mathcal{S}(n_1,\dots,n_k)$ is said to be a \emph{full affine
system of supports}.
\end{Def}

In the next Lemma we show that systems of supports are, in some
sense, closed under projections.

\begin{Lemma}\label{inductive} Let $k>1$, and let $(n_1,\dots, n_k)\in \N ^k$. Let $\Scal (n_1,\dots ,n_k)=(\Scal ;A_I, I\in \Scal)$ be
a system of supports. Fix $I\in \Scal \setminus \{\{1,\dots ,k\}\}$, and
let $p\colon (\No ^*)^k\to (\No ^*)^{\{1,\dots ,k\}\setminus I}$
denote the canonical projection.

If we define $\Scal _I=\{K\setminus I\mid K\in \Scal \mbox{ and }
I\subseteq K\}$ and for each $K\setminus I\in \Scal _I$ we take
$A_{I,K\setminus I}=A_K$, then:
\begin{itemize}
\item[(1)] $\Scal _I(p(n_1,\dots ,n_k))=(\Scal _I; A_{I,K\setminus
I},K\setminus I\in \Scal _I)$ is a system of supports of
$(\No^*)^{\{1,\dots ,k\}\setminus I}$.
\item[(2)] If $I\neq \emptyset$ then $\mid \Scal _I\mid <\mid \Scal\mid$.
\item[(3)] If $\Scal(n_1,\dots ,n_k)$ is a full affine system of supports then so is $\Scal _I(p(n_1,\dots
,n_k))$.
\end{itemize}
\end{Lemma}

\begin{Proof} It is routine to check that $\Scal _I(p(n_1,\dots
,n_k))$ satisfies the conditions of a system of supports. Statements
(2) and (3) are immediate from the definitions.
\end{Proof}

We note that the systems of supports is just an alternative way to
describe a class submonoids of $(\No^*)^k$,

\begin{Prop} \label{nonsense}
Fix $k \in \mathbb{N}$ and $(n_1,\dots,n_k) \in \mathbb{N}^k$. For
any $I\subseteq \{1,\dots ,k\}$, let $p_I \colon
(\mathbb{N}_0^*)^k \to (\mathbb{N}_0^*)^{\{1,\dots,k\}\setminus
I}$ denote  the canonical projection. Let
$\mathcal{S}(n_1,\dots,n_k)$ be a system of supports. Consider the
subset $M(\mathcal{S})$ of $(\No ^*)^k$ defined by $x\in
M(\mathcal{S})$ if and only if $I = \infsupp (x) \in \cal S$ and
$p_I(x) \in A_I$. Then $M(\mathcal{S})$ is a submonoid of
$(\No^*)^k$ such that  $(n_1,\dots,n_k)\in \N^k\cap M(\mathcal{S})$ and satisfying
the properties:
\begin{itemize}

\item[(M1)] if $I\subseteq \{1,\dots ,k\}$ is an infinite support of
some $x\in M(\mathcal{S})$ then the element $x^*$ determined by
$\supp (x^*) =\infsupp (x^*)=\infsupp (x)$ belongs to
$M(\mathcal{S})$.
\item[(M2)] If $x\in M(\mathcal{S})$ then $\infty\cdot x\in M(\Scal
)$.
\end{itemize}

Any submonoid $M$ of $(\No^*)^k$ with $M^k\cap \N ^k\neq
\emptyset$ and satisfying (M1) and (M2) is of the form
$M=M(\mathcal{S})$ for some system of supports $\mathcal{S}$.
Moreover, in this situation,   for any $I\in \mathcal{S}$,
$A_I=p_I(M)\cap \No ^{\{1,\dots,k\}\setminus I}$.
\end{Prop}

\begin{Proof}
Since $A_\emptyset \subseteq M(\mathcal{S})$, $0\in M(\mathcal{S})$
and $(n_1,\dots,n_k)\in M(\mathcal{S})$. To see that
$M(\mathcal{S})$ is a monoid it remains to see that it is closed
under addition.

Let $x,y\in M(\mathcal{S})$.  Set $I=\infsupp (x)$, $J=\infsupp
(y)$ and $K=I\cup J$. Let  $x^*\in (\No^*)^k$ be such that $\supp
(x^*) =\infsupp (x^*)=K$.  Notice that $x^*\in M(\mathcal{S})$
because $0\in A_K$, and that $x+x^*$ and $y+x^*\in M(\mathcal{S})$
because of condition (iv) in the definition of system of supports.
Then $x+y=(x+x^*)+(y+x^*)$ and, since $\infsupp (x+x^*)=\infsupp
(y+x^*)=K$, we deduce that $x+y\in M(\mathcal{S})$.

Note that, since the $A_I$ are monoids,  $M(\mathcal{S})$
satisfies (M1). Property (M2) follows from  condition (iii) in the
definition of system of supports.

To prove the converse, let $M$ be a submonoid of $(\No^*)^k$
satisfying (M1) and (M2) and such  that $(n_1,\dots ,n_k) \in
M\cap \N ^k$. Set
$$\mathcal{S}=\{I\mid \mbox{ there exists } x\in M \mbox{ such
that }\infsupp (x)=I\}$$ as the collection of subsets of $\{1,\dots
,k\}$. Moreover,  set $A_I=p_I(M)\cap \No ^{\{1,\dots,k\}\setminus
I}$ for any $I\in \mathcal{S}$. It is easy to check that  the
properties of $M$ ensure that $\mathcal{S}(n_1,\dots ,n_k)$ is a
system of supports such that $M(\mathcal{S})=M$.
\end{Proof}

\begin{Rem}\label{description} With the notation as in
Proposition~\ref{nonsense} and  Lemma~\ref{inductive}, assume that
$|\Scal|>2$ and let $\mathcal{T}=\{I_1,\dots ,I_\ell \}$ be the set
of minimal elements in $\Scal \setminus \{\emptyset\}$. For each $i \in
\{1,\dots ,\ell\}$, let $S_{I_i}$ be the system of supports given by
Lemma~\ref{inductive} and let $M(S_{I_i})$ be the associated monoid.
Then
$$M(\Scal)=A_\emptyset\bigcup\left(\cup _{i=1}^\ell
M'(S_{I_i})\right)=M_0\bigcup\left(\cup _{i=1}^\ell
M'(S_{I_i})\right)$$ where $M_0=A_\emptyset + \{\infty\cdot x\mid
x\in A_\emptyset \}$ and $$M'(S_{I_i})=\{x\in (\No^*)^k\mid
p_{I_i}(x)\in M(S_{I_i})\mbox{ and } \infsupp (x)\supseteq
I_i\}.$$
\end{Rem}

Now we give a couple of crucial examples of monoids given by a full
affine system of monoids.

\begin{Ex}\label{semilocalsystem} Let $k\ge 1$. Let $R$ be a noetherian semilocal ring with an onto
ring homomorphism $\varphi \colon R\to M_{n_1}(D_1)\times \cdots
\times M_{n_k}(D_k)$ with kernel $J(R)$ and where $D_1,\dots ,D_k$
are division rings. Let   $M=\dim _{\varphi}V^*(R)\subseteq
(\No^*)^k$. Then $M$ is a submonoid of $(\No^*)^k$ given by a full
affine system of supports.
\end{Ex}

\begin{Proof} Notice that  $(n_1,\dots,n_k)\in M$. By Theorem~\ref{behaviour}(1), $M$ satisfies
condition $(M1)$ in Proposition~\ref{nonsense}.

If $P$ is a countably generated projective right $R$-module then
$\dim _{\varphi}(\langle P^{(\omega)}\rangle)=\infty \cdot \dim
_{\varphi}(\langle P\rangle)$. Hence $M$ also satisfies condition
$(M2)$ in Proposition \ref{nonsense} and we can conclude that $M$
is a  submonoid of $(\No^*)^k$ given by a system of supports.

Let $x\in M$, and let $P$ be a countably generated projective right
such that $\dim _{\varphi}(\langle P\rangle)=x$. Let $I=\infsupp\,
(x)$. By Theorem~\ref{behaviour}(1), there exists   a
countably generated projective right $R$ module $P'$ such that
\[\supp (\dim
_{\varphi}(\langle P'\rangle))=\infsupp (\dim _{\varphi}(\langle
P'\rangle))=I.\] Let $J$ be the trace ideal of $P'$. Then $R/J$ is a
semilocal ring with Jacobson radical $J+J(R)/J$ and, by Lemma
\ref{simplefactors},
 $\varphi$ induces an onto ring homomorphism
$\overline{\varphi}\colon R/J\to \prod _{i\in \{1,\dots ,k
\}\setminus I}M_{n_i}(D_i)$ with kernel $J(R/J)$. Moreover, by
Theorem~\ref{behaviour}(2),
\[\dim
_{\overline{\varphi}}(\langle P/PJ\rangle)=p_{I}(x)\] where $P/PJ$
is a finitely generated projective right $R/J$-module. This shows
that $A_I=p_I(M)\cap \No^k\subseteq
\dim_{\overline{\varphi}}V(R/J)$. We claim that
$\dim_{\overline{\varphi}}V(R/J)=A_I$. Equivalently, for any
finitely generated projective right $R/J$-module $\overline P$ there
exists a countably generated projective right $R/J$-module $P_1$
such that $P_1/P_1J\cong \overline P$. Therefore the claim follows
from Proposition \ref{symtraces}(iii).

By Corollary~\ref{modjr}, $A_I=\dim_{\overline{\varphi}}V(R/J)$ is a
full affine submonoid of $\No ^{\{1,\dots ,k\}\setminus I}$.
Therefore the monoid $M$ is given by a full affine system of
supports as we wanted to show.
\end{Proof}

Next example is a consequence of Example~\ref{semilocalsystem} and
Theorem~\ref{pbrealization}. We prefer to give a proof just in the
monoid context.

\begin{Ex}\label{systemsex} Let $D\in M_{n\times k}(\No)$ and $E_1,$ $E_2\in M_{\ell \times
k}(\No)$. Let $m_1,\dots ,m_n\in \N$ be such that $m_i\ge 2$ for any
$i\in \{1,\dots ,n\}$.

Let $M\subseteq (\No^*)^k$ be the set of solution in $(\No^*)^k$ of the
system
\[D\left(\begin{array}{c}t_1\\\vdots \\ t_k\end{array}\right)\in \left(\begin{array}{c}m_1\No ^*
\\\vdots \\ m_n\No^*\end{array}
\right)\qquad \mbox{and} \qquad
E_1\left(\begin{array}{c}t_1\\\vdots \\
t_k\end{array}\right)=E_2\left(\begin{array}{c}t_1\\\vdots \\
t_k\end{array}\right).\] Assume also that $(n_1,\dots ,n_k) \in
M\cap \N ^k$. Then  there exists a full affine system of supports
$\mathcal{S}(n_1,\dots ,n_k)$ such that $M=M(\mathcal{S})$.
\end{Ex}

\begin{Proof}
By Lemma~\ref{infsupp}, $M$ satisfies  conditions (M1) and (M2) in
Proposition \ref{nonsense}.  Hence, by Proposition \ref{nonsense}
and following the notation there, $M$ is given by a system of
supports in which
$$\mathcal{S}=\{I\mid \mbox{ there exists } x\in M \mbox{ such that
}\infsupp (x)=I\}$$ as a collection of subsets of $\{1,\dots ,k\}$
and, for any $I\in \Scal$, $A_I=p_I(M)\cap \No
^{\{1,\dots,k\}\setminus I}$. We want to show that, for each $I \in
\mathcal{S}$, the monoid $A_I$ is full affine in $\No
^{\{1,\dots,k\} \setminus I}$. To this aim we prove that $A_I$ is
the set of solution in $\No ^{\{1,\dots,k\} \setminus I}$ of a
certain subsystem of the initial one.

Fix $I\in \mathcal{S}$. Let $D_I$ be the matrix with entries in
$\No$ obtained from $D=(d_{ij})$ by first deleting  the rows $i$
such that there exists $j\in I$ with $d_{ij}\neq 0$,  and then
deleting in the remaining matrix the $j$-th column for any $j\in I$.
Let us denote $K$ the subset of $\{1,\dots,n\}$ indicating which
rows of $D$ were deleted and let $p_K \colon {\N_0^*}^n \to
{\N_0^*}^{\{1,\dots,n\}\setminus K}$ be the canonical projection.

As $I \in \mathcal{S}$, for all $i \in \{1,\dots,\ell\}$ is either
$e_{i,j}^1 = 0 = e_{i,j}^2$ for all $j \in I$ or there are $j_1,j_2 \in I$
such that $e_{i,j_1}^1\neq 0$ and $e_{i,j_2}^2 \neq 0$.
Let $E_1^I$ and $E_2^I$ be the matrices with entries in $\No$
obtained from $E_1=(e^1_{ij})$ and $E_2=(e^2_{ij})$ by first
deleting the rows $i$ such that there exists $j\in I$ satisfying
that $e^1_{ij}$ is different from zero; after
we also delete to each of the remaining matrices the $j$-th column
for any $j\in I$.

Then the monoid $A_I$ is the set of solutions in
$\No^{\{1,\dots,k\} \setminus I}$ of the system
\[D_I\cdot p_I\left(\begin{array}{c}t_1\\\vdots \\
t_k\end{array}\right)\in p_K
\left(\begin{array}{c}m_1\No^*\\\vdots \\ m_n\No^*\end{array}
\right)\qquad \mbox{and} \qquad
E^I_1\cdot p_I\left(\begin{array}{c}t_1\\\vdots \\
t_k\end{array}\right)=E^I_2\cdot p_I\left(\begin{array}{c}t_1\\\vdots \\
t_k\end{array}\right).\]
\end{Proof}

\begin{Th}\label{affinesystem} Let $k\in \N$,  and let $(n_1,\dots
,n_k)\in \N ^k$.  Let $M$ be a submonoid of $(\No^*)^k$ such that
$(n_1,\dots ,n_k)\in M$. Then the following statements are
equivalent
\begin{enumerate}
\item[(i)] There exists a full affine system of supports $\mathcal{S}(n_1,\dots
,n_k)$ such that $M=M(\mathcal{S})$.
\item[(ii)] The submonoid $M$ is defined by a system of equations.
\end{enumerate}
\end{Th}

To easy the proof of the theorem we first prove an auxiliary result.

\begin{Lemma}
\label{hide} Let $(A,(n_1,\dots ,n_k))\subseteq (\No^*)^k$ be a
monoid defined by a system of equations $\mathcal{E}_A$. Let $I$
be a proper nonempty subset of $\{1,\dots ,k\}$ and set
$I_c=\{1,\dots ,k\}\setminus I$. Let $p\colon (\No^*)^k\to
(\No^*)^{I_c}$ denote the canonical projection. Assume also
$B\subseteq (\No ^*)^{I_c}$ is a monoid defined by a system of
equations $\mathcal{E}_B$ and such that $p(n_1,\dots ,n_k)\in B$.

Then the set $A'\subseteq (\No^*)^k$ defined by $x=(x_1,\dots, x_k
)\in A'$ if and only if either
\begin{enumerate}
\item[(1)] $x\in A$ and $p(x)\in B$
\end{enumerate}
or \begin{enumerate} \item[(2)] $x_i=\infty$ for any $i\in I$ and
$p(x)\in B$
\end{enumerate}
is a monoid defined by a system of equations. Moreover,
$(n_1,\dots ,n_k)\in A'$.
\end{Lemma}

\begin{Proof} Let \[D\left(\begin{array}{c}t_1\\\vdots \\ t_k\end{array}\right)\in \left(\begin{array}{c}m_1\No ^*
\\\vdots \\ m_n\No^*\end{array}
\right)\qquad \mbox{and} \qquad
E_1\left(\begin{array}{c}t_1\\\vdots \\
t_k\end{array}\right)=E_2\left(\begin{array}{c}t_1\\\vdots \\
t_k\end{array}\right)\] be a system of equations defining $A$. Fix
$i\in I$. Let $\mathcal{E}_A'$ be the system of equations
\[D\left(\begin{array}{c}t_1+n_1t_i\\\vdots \\ t_k+n_kt_i\end{array}\right)\in
\left(\begin{array}{c}m_1\No ^*
\\\vdots \\ m_n\No^*\end{array}
\right)\qquad \mbox{and} \qquad
E_1\left(\begin{array}{c}t_1+n_1t_i\\\vdots \\
t_k+n_kt_i\end{array}\right)=E_2\left(\begin{array}{c}t_1+n_1t_i\\\vdots \\
t_k+n_kt_i\end{array}\right).\] Notice that if $x\in \No^k$  then
it is a solution of $\mathcal{E}_A'$ if and only if $x\in A$;
while any $x\in (\No^*)^k$ such that $i\in \infsupp(x)$ is a
solution of $\mathcal{E}_A'$.

Let $A_i\subseteq (\No^*)^k$ be the solutions of the system
$\mathcal{E}_i=\mathcal{E}_A'\cup \mathcal{E}_B$, where
$\mathcal{E}_B$ is the trivial extension of the system defining
$B$. Then $x\in \No^k$   belongs to $A_i$ if and only if $x\in A$
and $p(x)\in B$. In particular, $(n_1,\dots,n_k)\in A_i$. While if
$x\in (\No^*)^k$ is such that $i\in \infsupp(x)$ then $x\in A_i$
if and only if $p(x)\in B$.

Now the submonoid $A'=\cap _{i\in I}A_i\subseteq (\No^*)^k$
defined as the solutions of the system $\cup _{i\in
I}\mathcal{E}_i$, satisfies the properties required in the
conclusion of the statement.
\end{Proof}

\begin{ProofTh} In view of Example~\ref{systemsex} we only need to prove
that (i) implies (ii).

Let $k \in \mathbb{N}$, and let $\mathcal{S}(n_1,\dots ,n_k)$ be a
full affine system of supports. We proceed by induction on $|\cal
S|$. If $|{\cal S}| = 2$ then the only  sets in $\mathcal{S}$ are
$\emptyset$ and $\{1,\dots, k\}$. So that the only nontrivial full
affine semigroup is $A_\emptyset$ and, for any $x\in
A_\emptyset\setminus \{0\}$, $\supp (x)=\{1,\dots, k\}$. By
\cite[Example~2.5]{FH2}, this implies that there exists
$y=(y_1,\dots ,y_k)\in A_\emptyset$ such that $A_\emptyset=y\No$.
Therefore $M(\mathcal{S})=A\cup \{\infty\cdot y\}$. By
Proposition~\ref{fas}, $M(\mathcal{S})$ is defined by a system of
equations.

Now assume that $|{\cal S}| > 2$ and that the statement  is true
for full affine systems of supports such that the set of supports
has smaller cardinality. Let $\cal T \subseteq \cal S$ be the set
of all minimal elements of $\cal S\setminus \emptyset $. Note
that, since $|{\cal S}| > 2$, no element in $\mathcal{T}$ is equal
to $\{1,\dots ,k\}$. For any $I \in \cal T$ we construct the full
affine system of supports $\Scal _I$ given by
Lemma~\ref{inductive}. As, for each $I \in \cal T$, $|\Scal
_I|<|{\cal S}|$, we know that the monoid $M_I(\Scal _I)$ is given
by a system of equations. Moreover, by Proposition~\ref{fas},
$M_0=A_\emptyset + \{\infty\cdot x\mid x\in A_\emptyset\}$ is a
submonoid of $(\No^*)^k$ defined by a system of equations.

Let ${\cal T} = \{I_1,\dots,I_\ell\}$. We complete $M_0$ to a chain
$M_0\subseteq M_1\subseteq \cdots \subseteq M_\ell$ of submonoids of
$(\No^*)^k$ given by a system of equations, inductively, in the
following way: If $i<\ell$ is such that $M_i$ is constructed then
$M_{i+1}$ is the monoid given by applying Lemma~\ref{hide} to
$A=M_i$, $I=I_{i+1}$ and $B=M_{I_{i+1}}(\Scal _{I_{i+1}})$. Notice
that, following the notation of  Lemma~\ref{hide} and by the
definition of a system of supports, $p_{I_{i+1}}(M_i)\subseteq
M_{I_{i+1}}(\Scal _{I_{i+1}})$, therefore $M_{i+1}=M_{i}\bigcup
M'_{I_{i+1}}$ where
$$M'_{I_{i+1}}=\{x\in (\No^*)^k\mid p_{I_{i+1}}(x)\in M_{I_{i+1}}(\Scal _{I_{i+1}})
\mbox{ and } \infsupp (x)\supseteq I_{i+1}\}.$$ Therefore
$$M_\ell =M_0\cup M'_{I_1}\cup \cdots \cup M'_{I_\ell}=M(\Scal)$$
by Remark~\ref{description}. This allows us to conclude that
$M(\Scal)$ is  a monoid given by a system of equations.
\end{ProofTh}

Now Theorem \ref{main} follows by just patching together our
previous results.

\bigskip

\begin{Proofmain}
$(1)\Rightarrow (2)$ follows from Theorem \ref{pbrealization}. It
is clear that $(2)\Rightarrow (3)$.

Finally, assume $(3)$. By Example \ref{semilocalsystem}, $M$ is
given by a full affine system of supports. By Theorem
\ref{affinesystem}, $M$ is defined by a system of equations, and
$(1)$ follows.
\end{Proofmain}

\section{Some consequences and  examples} \label{final}

In order to be able to construct further examples, first we single
out the monoids corresponding to semilocal rings such that any
projective right $R$-module is a direct sum of finitely generated
modules. They are precisely the ones arising in
Proposition~\ref{fas}.

\begin{Cor} \label{allfg} Let $k\in \N$. Let $M$ be a submonoid of $(\No^*)^k$
containing  $(n_1,\dots ,n_k)\in \N ^k$. Let $A=M\cap \No^k$. Then
the following statements are equivalent:
\begin{itemize}
\item[(i)] $A$ is a full affine submonoid of $\No^k$ and  $M=A+\{\infty\cdot a\mid a\in A\}$.
\item[(ii)] There exists a (noetherian) semilocal
ring $R$ such that all projective right $R$ modules are direct sum
of finitely generated modules,  and an onto morphism $\varphi \colon
R\to M_{n_1}(D_1)\times \dots \times M_{n_k}(D_k)$ with
$\mathrm{Ker}\, \varphi =J(R)$, where $D_1,\dots,D_k$ are suitable
division rings,  satisfying that $\dim _\varphi V^*(R)=M$.
\end{itemize}
\end{Cor}

\begin{Proof} Assume $(i)$. By Proposition \ref{fas}, $M$ is given by a system
of equations. By Theorem \ref{main}, $(ii)$ holds.

Assume $(ii)$. So that $R$ is a, not necessarily noetherian,
semilocal ring such that all projective right modules are direct sum
of finitely generated ones. By Corollary~\ref{modjr}, $A$ is a full
submonoid of $\No^k$. It is clear that $A+\{\infty\cdot a\mid a\in
A\}\subseteq M$. Let $P_1,\dots ,P_s$ be a set of representatives of
the indecomposable (hence finitely generated) projective right
modules. For $i=1,\dots ,s$, let $a_i=\dim_\varphi (\langle
P_i\rangle)$. As any projective module is a direct sum of finitely
generated projective modules, any $x\in M$ satisfies that $x=\alpha
_1a_1+\cdots +\alpha _sa_s$ for some $\alpha _i\in \No^*$, hence
$x\in A+\{\infty\cdot a\mid a\in A\}$. This shows that $(i)$ holds.
\end{Proof}

In Examples~\ref{noniso} we shall see that the property that all
projective modules  are direct sum of finitely generated ones is not
left right symmetric. By Theorem~\ref{behaviour} it is  left-right
symmetric in the noetherian semilocal case.

We stress the fact that, in view of Theorem~\ref{main}, noetherian
semilocal rings can have a rich supply of infinitely generated
projective modules that are not direct sum of finitely generated
ones.

It is quite  an interesting question to determine the structure of
$V^*(R)$ for a general semilocal ring. But right now it seems to be
too challenging  even for   semilocal rings  $R$ such that
$R/J(R)\cong D_1\times D_2$ for   $D_1$, $D_2$ division rings. Now
we provide some examples of such rings to illustrate
Theorem~\ref{main} and the difficulties that  appear in the general
case. We first observe that with such rings, since $k=2$ and $(1,1)$
is the order unit of $V(R)$, to have some room for interesting
behavior all finitely generated projective modules must be free.

\begin{Lemma}\label{free} Let $R$ be a semilocal ring such that $R/J(R)\cong D_1\times D_2$ for
suitable division rings $D_1$ and $D_2$. Fix $\varphi \colon R\to
D_1\times D_2$ an onto ring homomorphism such that $\mathrm{Ker}\,
\varphi =J(R)$. If $R$ has nonfree finitely generated projective
right (or left) modules then there exists $n\in \N$ such that $\dim
_\varphi V(R)$ is the submonoid of $\No ^2$ generated by $(1,1)$,
$(n,0)$ and $(0,n)$. In this case,
\[\dim _\varphi V^*(R)=(1,1)\No^*+(n,0)\No^*+(0,n)\No^*=\{(x,y)\in \No ^*\mid x+(n-1)y\in n\No ^*\}.\]
Therefore, all projective modules are direct sum of finitely
generated projective modules.
\end{Lemma}

\begin{Proof} Note that $\dim _\varphi (\langle R\rangle )=(1,1)$. So that $(1,1)\in A = \dim _\varphi
V(R)$.

Let $P$ be a nonfree finitely generated projective right $R$-module,
and let $\dim _\varphi (\langle P \rangle)=(x,y)$. As $P$ is not
free, either $x>y$ or $x<y$. Assume $x>y$, the other case is done in
a symmetric way. Then \[(x,y)=(x-y,0)+y(1,1)\in A\qquad (*)\] since,
by Corollary~\ref{modjr}, $A$ is a full affine submonoid of $\No ^2$
we deduce that $(x-y,0)\in A$ and also that
$(0,x-y)=(x-y)(1,1)-(x-y,0)\in A$. Choose $n\in \N$ minimal with
respect to the property $(n,0)\in A$, and note that then also
$(0,n)\in A$. We claim that
\[A=(1,1)\No+(n,0)\No +(0,n)\No.\]
We only need to prove that if $(x,y)\in A$ then it can be written as
a linear combination, with coefficients in $\No$ of $(1,1)$, $(n,0)$
and $(0,n)$. In view of the previous argument, it suffices to show
that if $(x,0)\in A$ then $(x,0)\in (n,0)\No$. By the division
algorithm $(x,0)=(n,0)q+(r,0)$ with $q\in \No$ and $0\le r<n$. As
$A$ is a full affine submonoid of $\No^2$ we deduce that $(r,0)\in
A$. By the choosing of $n$ we can deduce that $r=0$, as desired.

Let $P_1$ be a finitely generated right $R$-module such that $\dim
_\varphi (\langle P_1\rangle)=(n,0)$, and let $P_2$ be a finitely
generated right $R$-module such that $\dim _\varphi (\langle
P_2\rangle)=(0,n)$.

Let $Q$ be a countably generated projective right $R$-module that is
not finitely generated.  Let $\dim _\varphi (\langle Q\rangle
)=(x,y)\in \No ^*$. We want to show that
\[(x,y)\in (1,1)\No^*+(n,0)\No^*+(0,n)\No^*\]

If $x=y$ then $(x,y)=x(1,1)$ and, by Theorem~\ref{Pavel}(ii), $Q$ is
free. If $x>y$ then $y\in \No$ and $(x,y)=(x-y,0)+y(1,1)$, combining
Theorem~\ref{Pavel}(ii) with Lemma~\ref{full} we deduce that
$Q=yR\oplus Q'$ with $Q'$ such that $\dim _\varphi (\langle
Q'\rangle )=(z,0)$ where $z=x-y$. If $z<\infty$ then, by
Theorem~\ref{Pavel}(ii), $nQ'\cong zP_1$ hence $Q'$, and $Q$, are
finitely generated. So that, $z=\infty$ and then, by
Theorem~\ref{Pavel}(ii), $Q'\cong P_1^{(\omega)}$. Hence
$(x,y)=\infty \cdot (n,0)+y(1,1)$. The case $x<y$ is done in a
symmetric way.

It is not difficult to check that the elements of $\dim _\varphi
V^*(R)$ are the solutions in $\No^*$ of $x+(n-1)y\in n\No ^*$.
\end{Proof}

Now we  list all the possibilities for the monoid $V^*(R)$ viewed as
a submonoid of $V^*(R/J(R))$ when $R$ is a noetherian ring such that
$R/J(R)\cong D_1\times D_2$, for $D_1$ and $D_2$ division rings, and
all finitely generated projective modules are free.

We recall that a module is superdecomposable if it has no
indecomposable direct summand. By Theorem~\ref{main} and
Lemma~\ref{full}, in our context superdecomposable modules are
relatively frequent as they correspond to the elements $x$ in the
monoid $M\subseteq (\No^*)^k$ such that for any $y\in M\cap \No^k$
$\supp\, (y)\nsubseteqq \supp\, (x)$.

\begin{Ex} \label{nk2} Let $R$ be a semilocal noetherian ring such that there exists   $\varphi\colon R\to D_1\times D_2$, an
onto ring morphism with $\mathrm{Ker}\, \varphi=J(R)$, where  $D_1$
and $D_2$ are division rings. Assume that  all finitely generated
projective right $R$-modules are free. Hence $\dim _\varphi V(R)=
(1,1)\No$, and its order unit is $(1,1)$. Then there are the
following possibilities for $\dim_\varphi V^*(R)$:
\begin{itemize}
\item[(0)] All projective modules are free, so that
$M_0=\dim _\varphi V^*(R)=(1,1)\No^*$. Note that $M_0$ is the set of
solutions   $(x,y)\in \No^*$ of the equation $x=y$.
\item[(1)] $M_1=\dim _\varphi V^*(R)=(1,1)\No^*+ (0,\infty)\No^*$. So that, $M_1$ is the set of solutions $(x,y)\in \N^*$
of the equation $x+y=2y$.

Note that for such an $R$ there exists a countably generated
superdecomposable projective right $R$-module $P$ such that $\dim
_\varphi (\langle P\rangle )=(0,\infty)$. Then any countably
generated projective right $R$ module $Q$ satisfies that there
exist $n\in \No ^*$ and $m\in \{0,1\}$ such that $Q=R^{(n)}\oplus
P^{(m)}$.

\item[(1')] $M_1'=\dim _\varphi V^*(R)=(1,1)\No^*+ (\infty,0)\No^*$. So that, $M_1'$ is the set of solutions $(x,y)\in \N^*$
of the equation $x+y=2x$.

\item[(2)] $M_2=\dim _\varphi V^*(R)=(1,1)\No^*+ (\infty,0)\No^*+(0,\infty)\No^*$.
So that, $M_2$ is the set of solutions $(x,y)\in \No ^*$ of the
equation $2x+y=x+2y$.

Note that for such an $R$ there exist two  countably generated
superdecomposable projective right $R$-modules $P_1$ and $P_2$ such
that $\dim _\varphi (\langle P_1\rangle )=(0,\infty)$ and $\dim
_\varphi (\langle P_2\rangle )=(\infty ,0)$. Any countably generated
projective right $R$ module $Q$ satisfies that there exist $n\in
\No$ and $m_1,m_2\in \{0,1\}$ such that $Q=R^{(n)}\oplus
P_1^{(m_1)}\oplus P_2^{(m_2)}$.
\end{itemize}
\end{Ex}

\begin{Proof} In view of Theorem~\ref{main}, Theorem~\ref{affinesystem} and Lemma~\ref{free} we must describe all the possibilities for full
affine systems of supports of $\{1,2\}$ such that $A_{\emptyset} =
(1,1) \mathbb{N}_0$. Since the set of supports of a system of
supports at least contains $\emptyset$ and $\{1,2\}$ there are just
four possibilities.

Since the image of the projections of
$A_\emptyset$ either on the first or onto the second component is
$\N_0$, all the monoids $A_I$ in the definition of system of supports
are determined by $A_\emptyset$.

Case $(0)$ is the one in which $M_0=A_\emptyset+\infty\cdot
A_\emptyset$. According to Corollary \ref{allfg}, in this case all
projective modules are direct sum of finitely generated
(indecomposable) modules.

In cases $(1)$ and $(1')$ there are $3$ different infinite supports for the
elements in the monoid, and in case $(2)$ there are $4$.
\end{Proof}

The monoid $V^*(R)$ for non noetherian  rings $R$ with exactly two
maximal right ideals can have a more complicated structure. One of
the reasons is that  it may happen that $V(R)\subsetneqq W(R)$, cf.
\ref{modjr}. In the next Theorem we recall the almost only existing
family of such examples until now.

Note that, by Proposition~\ref{wr}, for a semilocal ring with a
fixed onto map $\varphi \colon R\to M_{n_1}(D_1)\times \cdots \times
M_{n_k}(D_k)$ to a semisimple artinian ring with kernel $J(R)$
\[W(R)=\{\langle P\rangle \in V^*(R)\mid \langle P/PJ(R)\rangle \in
V(R/J(R))\}\quad \mbox{and}\quad \dim_\varphi
W(R)=\left(\dim_\varphi V^*(R)\right)\cap \No^k.\] We point out the
fact that we do not even know whether the monoid $W(R)$ must be
finitely generated.

\begin{Th}\label{GS} Let $F$ be any field. Then there exists a semilocal $F$-algebra $R$ with an onto ring morphism  $\varphi\colon R\to F\times
F$ with $\mathrm{Ker}\, \varphi =J(R)$ and such that all finitely generated projective
modules are free but
\[\begin{array}{rcl}\dim _\varphi
W(R_R)&=&\{ (x,y)\in \No\mid x\ge y\}=(1,1)\No+(1,0)\No\\\dim
_\varphi V^*(R_R)&=&\left(\dim _\varphi W(R_R)\right)\No^*
\end{array}\] and
\[\begin{array}{rcl}\dim _\varphi W({}_RR)&=&\{ (x,y)\in \No\mid y\ge x\}=(1,1)\No+(0,1)\No \\ \dim
_\varphi V^*({}_RR)&=&\left(\dim _\varphi
W({}_RR)\right)\No^*\end{array}\]
\end{Th}

\begin{Proof} One possibility for such an example is the one
constructed by Gerasimov and Sakhaev in \cite{GS}. The monoids
$V^*(R_R)$ and $V^*({}_RR)$ were computed in \cite{DPP}.

Another possibility for such an example is the endomorphism ring of
a nonstandard uniserial module. Such modules were constructed by
Puninski in \cite{Pun1}, and it was observed in \cite{FHS} that the
endomorphism ring of such modules have such pathological projective
modules. The monoids were computed in \cite[Theorem~4.7]{P1}.
\end{Proof}

   Next family of
examples shows that the combination of  Theorem~\ref{GS} with the
pullback constructions gives a good tool to construct further
examples.

The duality between finitely generated projective right modules and
finitely generated projective left modules does not extend to an
isomorphism between $W(R_R)$ and $W({}_RR)$. But, as it follows from
Proposition~\ref{wr}, in the  semilocal case, one determines the
other.

In the next Example we show that, in general, for semilocal rings
the monoids $V^*({}_RR)$ and $V^*(R_R)$ not only are not isomorphic
but also that one monoid does not determine the other. We do that by
giving an example of a semilocal ring $R$ such that all  projective
left modules are free but $V^*(R_R)$ is isomorphic to the monoid in
Example \ref{nk2}(1'); in particular, not all right projective
modules are direct sum of finitely generated ones.

\begin{Ex}\label{noniso} Let $F$ be any field. In all the statements $R$ denotes a semilocal F-algebra, and
  $\varphi\colon R\to F\times F$ denotes an onto ring homomorphism
such that $\mathrm{Ker}\, \varphi=J(R)$. Fix $n\in \N$. Then:
\begin{itemize}
\item[(i)] There exist $R$ and $\varphi$ such that
\[N_1=\dim_\varphi V^*(R_R)=(1,1)\No^*+(n,0)\No^*=\{(x,y)\in (\No^*)^2\mid x\ge y \mbox{ and }x+(n-1)y\in n\No ^*\}\]
\[N_2=\dim_\varphi V^*({}_RR)=(1,1)\No^*+(0,n)\No^*=\{(x,y)\in (\No^*)^2\mid x\le y \mbox{ and }x+(n-1)y\in n\No ^*\}\]
\item[(ii)] There exist $R$ and $\varphi$ such that
\[\dim_\varphi V^*(R_R)=N_1+(0,\infty)\No^*=\{(x,y)\in (\No^*)^2\mid 2x+y\ge 2y+x \mbox{ and }x+(n-1)y\in n\No ^*\}\]
\[\dim_\varphi V^*({}_RR)=N_2+(\infty ,0)\No^*=\{(x,y)\in (\No^*)^2\mid 2x+y\le 2y+x \mbox{ and }x+(n-1)y\in n\No ^*\}\]
\item[(iii)] There exists $R$ and $\varphi$ such that
$\dim_\varphi V^*(R_R)=(1,1)\No^*+(\infty,0)\No^*$ and $\dim_\varphi
V^*({}_RR)=(1,1)\No^*$. Therefore, all projective left $R$-modules
are direct sum of finitely generated modules but this is not true
for projective right $R$-modules. In particular, $V^*(R_R)$ and
$V^*({}_RR)$ are not isomorphic.
\end{itemize}
\end{Ex}

\begin{Proof} To construct $(i)$, let $R_1$ be an $F$ algebra with
an onto ring homomorphism  $j_1\colon R_1\to D_1\times D_2$ where
$\mathrm{Ker}\, j_1=J(R)$ and $D_1$, $D_2$  are division rings.
Moreover, assume that $\dim
_{j_1}V^*(R)=(1,1)\No^*+(n,0)\No^*+(0,n)\No^*=M$, cf.
Lemma~\ref{free}. The existence of $R_1$ and $j_1$ is a consequence
either of Theorem~\ref{main} or of \cite[Theorem~6.1]{FH}.

Let $R_2$ and $j_2\colon R_2\to F\times F$ be ring and the onto map,
respectively, given by Theorem~\ref{GS}. That is,
\[
M_1=\dim _{j_2} V^*({R_2}_{R_2})=\{ (x,y)\in \No^*\mid x\ge
y\}=(1,1)\No^*+(1,0)\No^*\] and
\[M_2=\dim _{j_2} V^*({}_{R_2}R_2)=\{ (x,y)\in \No^*\mid y\ge x\}=(1,1)\No^*+(0,1)\No^* .\]

As $F\times F\subseteq D_1\times D_1$, we extend $j_2\colon R_2 \to
D_1\times D_2$ and consider the pullback
$$\begin{CD}
R_1@> j_1 >>D_1\times D_2\\
@A i_1 AA @AA j_2 A\\
R@>> \varphi > R_2
\end{CD}$$

By Corollary~\ref{milnorsemilocal} and
Corollary~\ref{milnorintersection}, $R$ is a semilocal $F$-algebra,
$\varphi$ is an onto ring homomorphism with kernel $J(R)$ and
\[\dim_\varphi V^*(R_R)=M\cap M_1\qquad \qquad \dim_\varphi
V^*({}_RR)=M\cap M_2\] which are the monoids $N_1$ and $N_2$
respectively.

To construct example $(ii)$  set $R_1$ and $j_1\colon R_1\to F\times
F$ to be the ring and the onto ring homomorphism, respectively,
constructed in $(i)$. Let $j_2\colon F\times F\to M_3(F)\times
M_3(F)$ be the map defined by
\[j_2(x,y)=\left(\left(\begin{smallmatrix}x&0&0\\0&x&0\\0&0&y\end{smallmatrix}\right),
\left(\begin{smallmatrix}y&0&0\\0&x&0\\0&0&y\end{smallmatrix}\right)\right)\]
Note that $j_2$ induces the monoid homomorphism $f\colon
(\No^*)^2\to (\No^*)^2$ defined by $f(x,y)=(2x+y, 2y+x)$.

Let $R$ be the ring defined by the pullback
$$\begin{CD}
M_3(R_1)@> M_3(j_1) >>M_3(F)\times M_3(F)\\
@A i_1 AA @AA j_2 A\\
R@>> \varphi >  F\times F
\end{CD}$$

By Corollary~\ref{milnorsemilocal},  $R$ is a semilocal $F$-algebra,
$\varphi$ is an onto ring homomorphism with kernel $J(R)$, and
$\dim_\varphi V^*(R_R)$ and $\dim_\varphi V^*({}_RR)$ are as claimed in $(ii)$.

To construct example $(iii)$, let $R_1$ and $j_1\colon R_1\to
F\times F$ be ring and the onto map, respectively, given by
Theorem~\ref{GS}. Let $R_2$ be a  semilocal ring with an onto ring
homomorphism $j_2\colon R_2\to F\times F$, with $\mathrm{Ker}\,
j_2=J(R_2)$, satisfying that $\dim _{j_2} V^*((R_2)_{R_2})=\dim
_{j_2} V^*({}_{R_2}R_2)=(1,1)\No^*+(\infty,0)\No^*$. Note that such
example exists by Example~\ref{nk2}(1).

Let $R$ be the ring defined by the pullback
$$\begin{CD}
R_1@> j_1 >>F\times F\\
@A i_1 AA @AA j_2 A\\
R@>> \varphi >  R_2
\end{CD}$$

By Corollary~\ref{milnorintersection},  $R$ is a semilocal
$F$-algebra, $j_2\varphi$ is an onto ring homomorphism with kernel
$J(R)$ and
\[\dim_{j_2\varphi} V^*(R_R)=\left((1,1)\No^*+(1,0)\No^*\right)\cap \left((1,1)\No^*+(\infty,0)\No^*\right)= (1,1)\No^*+
(\infty,0)\No^*\]
\[\dim_{j_2\varphi} V^*({}_RR)=\left((1,1)\No^*+(0,1)\No^*\right)\cap \left((1,1)\No^*+(\infty,0)\No^*\right)= (1,1)\No^*\]
Therefore all projective left $R$-modules are free. On the other
hand if $P$ is a countably generated projective right $R$-module
such that $\dim_{j_2\varphi}(\langle P\rangle)=(\infty,0)$ then $P$
is superdeccomposable, hence, it is not a direct sum of finitely
generated projective right $R$-modules. In order to see that $V^*(R_R)$ and $V^*(_RR)$ are not isomorphic, it is
enough to count idempotents in these monoids.
\end{Proof}

\begin{Rem}{\rm
Example \ref{noniso}(ii) answers a problem mentioned in
\cite[page~3261]{DPP}, and Example \ref{noniso}(iii) answers a
problem in \cite[page~310]{FS}.}
\end{Rem}

It follows from Proposition~\ref{wr} that Corollary \ref{weakdiv}
fails for general semilocal rings, as the order induced on $W(R)$,
viewed as a submonoid of  $V^*(R/J(R))$, is not the algebraic order.
Next example shows that this can fail also for the subsemigroup
$W(R)\setminus V(R)$.

\begin{Ex}
\label{nodivisibility} For any field $F$, there exists a semilocal
$F$-algebra $R$ and an onto morphism $\varphi\colon R\to F\times
F\times F$, with $\mathrm{Ker}\, \varphi =J(R)$, satisfying that
\[M_1=\dim_\varphi V^*(R_R)=\{(x,y,z)\in (\No^*)^3\mid x\ge y\mbox{ and  } y\ge z\}\]
\[M_2=\dim_\varphi V^*({}_RR)=\{(x,y,z)\in (\No^*)^3\mid y\ge x\mbox{ and  } z\ge y\}\]
In particular, $a=(1,1,0)$ and $b=(1,0,0)\in M_1$, clearly $b<a$,
but there is no element $b'\in M_1$ such that $b+b'=a$.
\end{Ex}

\begin{Proof} Let $T$ be the $F$-algebra given by Theorem \ref{GS}. Set $R_1=T\times T$. Let
$j_1\colon R_1\to F^4$ be an onto ring homomorphism with kernel
$J(R_1)=J(T)\times J(T)$. Let $j_2\colon F\times F\times F\to F^4$
be the morphism defined by $j_2(x,y,z)=(x,y,y,z)$.

Note that $j_2$ induces the morphism of monoids $f\colon
(\No^*)^3\to (\No^*)^4$ also defined by $f(x,y,z)=(x,y,y,z).$

Let $R$ be the ring defined by the pullback diagram
$$\begin{CD}
R_1@> j_1 >>F^4\\
@A i_1 AA @AA j_2 A\\
R@>> \varphi > F\times F\times F
\end{CD}$$

By Corollary~\ref{milnorsemilocal}, $R$ is a  semilocal $F$-algebra
 and $\varphi$ is an onto ring homomorphism with kernel $J(R)$.
We have chosen $R_1$  such that $\dim _{j_1}
V^*((R_1)_{R_1})=\{(x,y,z,t)\in (\No^*)^4\mid x\ge y\mbox{ and }
z\ge t\}$. By Corollary~\ref{milnorsemilocal},
 $(x,y,z)\in \dim _{\varphi}V^*(R_R)$ if and only
if $x\ge y\ge z$  as desired.  Similarly, as $\dim
_{j_1}V^*({}_{R_1}{R_1})=\{(x,y,z,t)\in (\No^*)^4\mid x\le y\mbox{
and } z\le t\}$, $\dim_\varphi V^*({}_RR)=\{(x,y,z)\in (\No^*)^3\mid
x\le y\mbox{ and  } y\le z\}$. The rest of the statement is clear.
\end{Proof}

\renewcommand{\baselinestretch}{1}

\end{document}